\documentclass[11pt] {article}
  \usepackage{amsmath}
    \usepackage{amssymb}
    
     \usepackage{pdfsync}

\newtheorem{theorem}{Theorem}[section]
\newtheorem{lemma}{Lemma}[section]

\newtheorem{remark}{Remark}[section]

\newcommand{\eqnsection}{
   \renewcommand{\theequation}{\thesection.\arabic{equation}}
   \makeatletter
   \csname @addtoreset\endcsname{equation}{section} 
   \makeatother}

\def \ov{\overline}

\def \be{\begin{equation}}
\def \ee{\end{equation}}
\def \bt{\begin{theorem}} 
\def \et{\end{theorem}}
\def \bl{\begin{lemma}} 
\def \el{\end{lemma}}
\def \bea{\begin{eqnarray}}
\def \eea{\end{eqnarray}}
\def \bas{\begin{eqnarray*}}
\def \eas{\end{eqnarray*}}

\def \al{\alpha}
\def \bb{\beta}
\def \ga{\gamma}

\def \de{\delta}

\def \ep{\epsilon}

\def \la{\lambda}

\def \om{\omega}
\def \Om{\Omega}

\def \si{\sigma}

\def \ff{\infty}

\def \wt{\widetilde}

\def\stl{\stackrel{law}{=}}
\def\ds{\displaystyle}

\def \AA{{\cal A}}
\def \BB{{\cal B}}

\def \JJ{{\cal J}}

\def\b1{\mathbf 1}

\def \({\left(}
\def \){\right)}

\def \nn{\nonumber}
\def \Proof{\noindent{\bf Proof $\,$ }}

\def \bc{\begin{center} }
\def \ec{\end{center} }
\def \bs{\begin{slide} }
\def \es{\end{slide} }

\def\square{{\vcenter{\vbox{\hrule height.3pt
        \hbox{\vrule width.3pt height5pt \kern5pt
           \vrule width.3pt}
        \hrule height.3pt}}}}
\def\qed{{\hfill $\square$ \bigskip}}

\eqnsection

\begin{document}

\title{ Permanental sequences  related to
  a Markov chain  example of Kolmogorov }

  \author{  Michael B. Marcus\,\, \,\, Jay Rosen \thanks{Research of     Jay Rosen was partially supported by a grant from the Simons Foundation.   }}
\maketitle

 \begin{abstract}   Permanental sequences with non-symmetric kernels that are  
generalization of the potentials of a Markov chain with state space $\{0,1/2,  \ldots, \newline 1/n,\ldots\}$ and a single instantaneous state   that was introduced  by Kolmogorov, are studied.  Depending on a   parameter in the  kernels we obtain  an exact rate of divergence of the sequence at $0$, an exact local modulus of continuity of the sequence at $0$, or a precise bounded discontinuity for   the sequence at $0$.

The kernel of the permanental sequence is,  
 \begin{equation}
U(0,0)=2,\quad U(0,1/j)=1+g_{j},  \hspace{.1 in}U(1/i,0)=1+ f_{i},\hspace{.1 in}i,j=2,3,\ldots .\label{a2}\nonumber
\end{equation}  
\vspace{-.2 in}   \begin{equation}
U(1/i,1/j)=\la_{j}\de_{i,j}+1+f_{i}g_{j}, \hspace{.2 in}\nn i,j=2,3,\ldots\label{a1},
\end{equation}
where  $\{\lambda_{j}\}$, $\{f_{i}\}$   and  $\{g_{j}\}$ satisfy certain conditions. 

Let $T=\{0,1/2,1/3,\ldots\}$. For all   $\alpha>0$,  there exist $\alpha$-permanental sequences    $\{X_{(\alpha), s};s\in T\}$ with kernel $U$. 
If 
$
\lim_{n\to\infty} \lambda_n\log n
 	=\beta $,  $0\le \beta<\infty,$
  then for all $k\ge 1$,
\begin{equation}
\limsup_{n\to \infty}\frac{X_{(k/2), 1/n}- X_{(k/2),0}}{(\lambda_n\log n )^{1/2}} =\beta^{1/2}+2 X_{(k/2),0}^{1/2}\,\quad a.s.,\label{21.pp}\nonumber
\end{equation}   \vspace{-.15in} 
  \begin{equation}
  \liminf_{ n\to \infty}\frac{ X_{(k/2), 1/n}- X_{(k/2),0} }{ ( \lambda_n\log n )^{1/2}}=  \left\{\begin{array} {cc}  -  X_{(k/2),0}/\beta^{1/2}& \mbox{if}\hspace{.05 in}  X_{(k/2),0}<\beta\\
\beta^{1/2} -2   X_{(k /2),0}^{1/2}& \mbox{if}\hspace{.05 in} X_{(k/2),0}\ge \beta
\end{array}\right.   \,a.s.\label{222}\nn
\end{equation}

  \end{abstract}  
 
\maketitle

\section{Introduction}\label{sec-1}

We  are interested in $\al$-permanental processes that are positive infinitely divisible processes determined by an infinite matrix that is the potential density of a transient Markov chain. When the matrix is symmetric a 1/2-permanental process it is the square of a Gaussian process. Permanental processes are related by the Dynkin isomorphism theorem to the total accumulated local time of the chain when the potential density is symmetric, and by a generalization of the Dynkin theorem by Eisenbaum and Kaspi in the general case. They are also related to  chi square processes and loop soups.

In this paper we study permanental sequences with non-symmetric kernels
that are a generalization of the potentials of  a   Markov chain   with a single instantaneous state that was introduced by
Kolmogorov [4]. The connection is explained  in Section \ref{sec-9}. These permanental sequences are very interesting because, suitably normalized,
 they can have bounded random discontinuities at zero. 

 Permanental processes   provide a challenge to probabilists who are interested in sample path properties of stochastic processes to see what new ideas and techniques are needed to analyze them. We think that they should have applications in statistical modeling that use chi square processes.

  An $R_{+}^{n}$ valued   $\al$-permanental random variable $X=(X_{1},\ldots, X_{n})$ is  a random variable with Laplace transform 
\begin{equation}
   E\(e^{-\sum_{i=1}^{n}s_{i}X_{i}}\) 
 = \frac{1}{ |I+KS|^{ \al}}   \label{int.1} ,
 \end{equation}
for some $n\times n$ matrix $K$, diagonal matrix $S$ with entries $s_{i}$, $1\le i\le n$, and $\al>0$.  
  We refer to $K$ as a kernel of $X$.      

     An $\al$-permanental process $\{X_{t},t\in T \}$ is a stochastic process  which has finite dimensional  distributions  that are $\al$-permanental vectors. The permanental process is determined by   a kernel $K=\{K(s,t),s,t\in T \}$, with the property that for all $t_{1},\ldots,t_{n}$ in $T $,  $\{K(t_{i},t_{j}),i,j\in [1,n] \}$ determines the $\al$-permanental random variable $ (X_{t_{1}},\ldots, X_{t_{n}})$ by (\ref{int.1}).   In this paper we take   $T =\{0,1/2,1/3, \ldots,1/n,\ldots\}$ with the Euclidean topology. 
     
      We refer to an $\al$-permanental process on $T$ as an $\al$-permanental sequence, or simply as a permanental sequence.  
 Note that when (\ref{int.1})  holds  for a kernel  $K(s,t)$ for all $\al>0$, the family of permanental processes obtained are infinitely divisible.

 It is well known that when $K$ is symmetric it is the covariance of a Gaussian sequence $\{\xi_{s}\}$ and $\Xi:=\{\Xi_{s},s\in T\}\stl \{\xi^{2}_{s}/2,s\in T\}$ is a $1/2$-permanental process. Consequently  it is relatively easy to find sample path properties of $\Xi$ because we have many tools to analyze Gaussian sequences. Moreover,  we can 
often extend these results to $k/2$-permanental sequences, for all integers $k\ge 1$, by considering the sum of $k$ independent copies of $\Xi$.  Therefore the real challenge is to understand the behavior  of  $\al$-permanental sequences for which the kernel $K$ in (\ref{int.1}) can not be taken to be  symmetric.  
What this means   is explained in detail in Section \ref{sec-9}. See also \cite{EM, MRns}.
  
 \medskip The kernel  $U = \{U(s, t), s, t \in T\}$ of the permanental sequences we consider  are obtained
as follows:
Let   $G=\{G_{i,j};i,j\in 0,2,\ldots \}$ where,
\be
G_{i,j}=\la_{i}\de_{i,j}+1+f_{ i}g _{j} ,\label{1.2}
\ee
 where $\la_{0}=0$,   $f_{ 0}=g _{0}=1$ and   $0<f_{i}, g_{i}<1$,   $i\ge 2$.  Written out this looks like,
\be    G= \left (
\begin{array}{   cccc c}  2&1+ g _{2}&\dots &1+ g _{n}&\dots \\
  1+f_{ 2} &\la_{2}+1+f_{ 2}g _{2}&\dots &1+f_{ 2}g _{n}&\dots \\
\vdots&\vdots&\ddots&\vdots &\ddots \\
  1+f_{ n} & 1+f_{ n}g _{2}  &\dots &   \la_{n}+1+f_{ n}g _{n}&\dots\\
\vdots&\vdots&\ddots&\vdots &\ddots\end{array}\right ).\label{21.0}
     \ee 
   We assume that $ \la_{j}>0$,  for $j\ge 2$ and,    \begin{equation}
\lim_{n\to \ff} \la_{n}=0.  \label{21.1}
     \end{equation}
 	We also take $\{f_{j}\}$, $\{g_{j}\}$ and $\{\la_{j}\}$ to satisfy,
\begin{equation}
\sum_{j=2}^{\ff} \frac{\(1- f_{j} \)}{\la_{j}}<1\quad\mbox{and}\quad \sum_{j=2}^{\ff} \frac{\(1- g_{j} \) }{\la_{j}}<1.\label{21.2}
\end{equation}

  Let $U(s,t)$ be a kernel on $T=\{0,1/2,1/3,\ldots, 1/n,\ldots\}$     defined  by,  
\begin{equation}
U(0,0)=G_{0,0}=2,\label{21.1c}
\ee
\be\quad U(0,1/j)=G_{1,j}=1+g_{j},  \hspace{.2 in}U(1/i,0)=G_{i,1}=1+ f_{i},\qquad i,j\ge 2 \nn
\end{equation}  
 \begin{equation}
\  U(1/i,1/j)=G_{i,j}=\la_{i}\de_{i,j}+1+f_{i}g_{j} \label{21.1b},\qquad i,j\ge 2.\nn
\end{equation}   
Note  that $U(s,t)$ is continuous on $T$.

 \begin{theorem}\label{theo-1.1mm} For all $\al>0$ there exists a permanental   sequence  $X_{(\al)}=\{X_{(\al), s} ,s\in T\} $ with kernel $U(s,t)$. 
\end{theorem}

We study the limiting behavior of $X_{(\al)}$ at 0. The next result is obtained in our earlier papers as we point out on  page \pageref{page-theo-1.2}. Its proof is much simpler than the proofs   of  the other limit theorems. 

  \bt\label{theo-1.2mm}
If  
\begin{equation}
\lim_{n\to \ff}\la_{ n}\log n
= \ff\label{21.3},
\end{equation}
 then for all $ \al>0$,
 \be
\limsup_{n\to \ff}\frac{ X_{(\al), 1/n}}{  \la_{n}\log
n  }=1\qquad 
a.s.\label{21.4}
\ee  
\et

 In   the next theorem  we provide examples of permanental processes with an almost surely bounded random discontinuity at $0$. Note that we are restricted to $\al=k/2$ for $k$ an integer greater that or equal to 1. We need this restriction to find   lower bounds. In Theorem \ref{theo-kolexp2s} we give upper bounds for all $\al>0$.

\bt\label{theo-kolexp12} 
If,  
\begin{equation}
\lim_{n\to \ff}\la_{n}\log
n= \bb, \hspace{.2 in}0<\bb<\ff,\label{21.3bb}
\end{equation}
 then for all $k\ge 1$,
\begin{equation}
\limsup_{  n\to \ff }|X_{(k/2), 1/n}- X_{(k/2),0}|=\bb+2\bb ^{1/2}X_{(k/2),0}^{1/2}\,\,, \hspace{.1 in}a.s.\label{21.4mv}
\end{equation}
\begin{equation}
\limsup_{  n\to \ff }X_{(k/2), 1/n}- X_{(k/2),0}=\bb+2\bb ^{1/2}X_{(k/2),0}^{1/2}\,\,, \hspace{.1 in}a.s.\label{21.4m}
\end{equation}
and    \begin{equation}
\liminf_{  n\to \ff }X_{(k/2), 1/n}- X_{(k/2),0}=  \left\{\begin{array} {cc}  -  X_{(k/2),0}, & \mbox{if }\hspace{.2 in}  X_{(k/2),0}< \bb\\
\bb -2 \bb^{1/2} X_{(k /2),0}^{1/2},& \mbox{if }\hspace{.2 in} X_{(k/2),0}\geq\bb
\end{array}\right.   , \hspace{.1 in}a.s.\label{21.4ms}
\end{equation}
\et

The next limit theorem gives a  local modulus of continuity   at zero that is itself random. This shows in particular that the modulus of continuity is not a tail event.

 	\bt \label{theo-kolexp22}

If, 
\begin{equation}
\lim_{n\to\ff}\la_{n}\log
n =0,\label{21.3cw}
\end{equation}
 then for all $k\ge 1$,  
 \begin{equation}
\limsup_{  n\to \ff }{|X_{(k/2), 1/n}- X_{(k/2),0}| \over  ( \la_{n}\log n)^{1/2}}=2 X^{1/2}_{(k/2),0}\,\,, \hspace{.1 in}a.s.\label{21.4n}
\end{equation}
\begin{equation}
\limsup_{  n\to \ff }{X_{(k/2), 1/n}- X_{(k/2),0} \over  ( \la_{n}\log n)^{1/2}}=2 X^{1/2}_{(k/2),0}\,\,, \hspace{.1 in}a.s.\label{21.4nq}
\end{equation}
and
\begin{equation}
\liminf_{  n\to \ff }{X_{(k/2), 1/n}- X_{(k/2),0} \over  (\la_{n}\log n)^{1/2}} =-2 X^{1/2}_{(k/2),0}\,\,, \hspace{.1 in}a.s.\label{21.4nn}
\end{equation}
  \et

 The next theorem is simply a restatement of Theorems \ref{theo-kolexp12} and \ref{theo-kolexp22} in a more compact form.

\bt\label{theo-kolexp2q} 
If, 
\begin{equation}
\lim_{n\to \ff} \la_{n}\log
n  = \bb, \hspace{.2 in}0\le \bb<\ff,\label{21.3b}
\end{equation}
 then for all $k\ge 1$,
\bea
\limsup_{ n\to \ff }\frac{|X_{(k/2), 1/n}- X_{(k/2),0}|}{( \la_{n}\log
n )^{1/2}}&=&\bb^{1/2}+2 X_{(k/2),0}^{1/2}\,,\qquad a.s.\label{21.4mva}
 \\&&\nn\\
 \limsup_{ n\to \ff }\frac{ X_{(k/2), 1/n}- X_{(k/2),0} }{ ( \la_{n}\log
n )^{1/2}} &=&\bb^{1/2}+2 X_{(k/2),0}^{1/2} \,, \qquad a.s.\label{21.4ma}
\eea
and
  \begin{equation}
  \liminf_{ n\to \ff }\frac{ X_{(k/2), 1/n}- X_{(k/2),0} }{ ( \la_{n}\log
n )^{1/2}}=  \left\{\begin{array} {cc}  -  X_{(k/2),0}/\bb^{1/2}, & \mbox{if }\hspace{.2 in}  X_{(k/2),0}<\bb\\
\bb^{1/2} -2   X_{(k /2),0}^{1/2},& \mbox{if }\hspace{.2 in} X_{(k/2),0}\ge \bb
\end{array}\right.   , \hspace{.1 in}a.s.\label{21.4msa}
\end{equation}
\et

   When $\bb>0$ in (\ref{21.3b}) we can replace the denominator on the left-hand sides  of (\ref {21.4mva})--(\ref{21.4msa})  by $\bb^{1/2}$ and then multiply by $\bb^{1/2}$ to get (\ref{21.4mv})--(\ref{21.4ms}).
 
  When $\bb=0$ in (\ref{21.3b}) we simply  replace $\bb$ by 0  on the right-hand sides  of (\ref {21.4mva})--(\ref{21.4msa})     to get (\ref{21.4n})--(\ref{21.4nn}).

\medskip	
Theorem \ref{theo-kolexp2q} is the main result in this paper. We note that it is   fully equivalent to the following theorem,  which has  a much simpler form but does not  make evident   the interesting types of    limits that we see in Theorem \ref{theo-kolexp2q}.  

\begin{theorem} \label{theo-1.4}  If,
\begin{equation} 
\lim_{n\to \ff} \la_{n}\log
n  = \bb, \hspace{.2 in}0\le \bb<\ff,\label{firstx}
\end{equation}
 then for all $k\ge 1$,
 \be 
\limsup_{n\to \ff}\frac{ |X^{1/2}_{(k/2), 1/n}- X^{1/2}_{(k/2),0}| }{( \la_{n}\log
n )^{1/2}}= 1\,,\qquad a.s.\label{firstqq}
 \ee
\be 
\limsup_{n\to \ff}\frac{ X^{1/2}_{(k/2), 1/n}- X^{1/2}_{(k/2),0} }{( \la_{n}\log
n )^{1/2}}= 1\,,\qquad a.s.\label{firstq}
 \ee

 \begin{equation}
  \liminf_{ n\to \ff}  \frac{ X^{1/2}_{(k/2), 1/n}- X^{1/2}_{(k/2),0} }{( \la_{n}\log
n )^{1/2}}=  
 - ({X_{(k /2),0}^{1/2}}/{\bb^{1/2}} \wedge 1)\,,\qquad a.s.\label{secondq}
  \end{equation}

 \end{theorem}

   \medskip	The next theorem also is in our earlier work. It shows  that the upper bounds given in   (\ref{21.4m}), (\ref{21.4nq}) and   (\ref{21.4ma}) hold for all  permanental process   $X_{(\al)}=\{X_{(\al),s};s\in T\}$  with kernel $U(s,t)$. The innovation in this paper is to find the lower bounds, which generally is more difficult.

    \bt\label{theo-kolexp2s} 
The upper bounds in (\ref{21.4m}), (\ref{21.4nq}) and equivalently (\ref{21.4ma}) hold for all $\al>0$.
\et

 The kernel $U$ defined in   (\ref{21.1c})    is related to the  potential of a Markov chain introduced by Kolmogorov. In    Kolmogorov's   paper \cite{Kol},  published in 1951, he gives an example of a recurrent Markov chain on the integers with a single instantaneous state.  This is generalized in Reuter \cite{Reu}, published in 1969. 
 These Markov chains are very interesting but the  potentials we use to  define permanental processes must be the potentials of  transient chains. Therefore, we modify  Reuter's example. We explain  this at the  end of this paper, in Section \ref{sec-9}, since our motivation for the definition of $U$ is not used to obtain the results given above.
 
 \medskip  Theorems   \ref{theo-1.1mm},  \ref{theo-1.2mm} and \ref{theo-kolexp2s}  are based on our earlier work. Their proofs are given in Section \ref{sec-5}.    Theorem \ref{theo-kolexp2q}  is the main result in this paper. We have commented in this section that it easily implies Theorems \ref{theo-kolexp12} and \ref{theo-kolexp22}. It is proved in Section \ref{sec-3}. Several critical lemmas used in the proof of Theorem \ref{theo-kolexp2q} are given Section \ref{sec-2}.  Lemmas \ref{lem-2.2}    and \ref{cor-7.1}    are proved in Section \ref{sec-2}.      Lemmas \ref{lem-kmatrix},    \ref{lem-koldet} and \ref{lem-gold} are proved in Section \ref{sec-isimi}.      The proof that  Theorem \ref{theo-kolexp2q} and  Theorem \ref{theo-1.4}  are equivalent is given in Section \ref{sec-4}. 

  \medskip We thank Pat Fitzsimmons for helpful conversations.

  \section{Preliminaries  }	\label{sec-2}
  Let  $C=\{ c_{ i,j}\}_{ 1\leq i,j\leq n}$ be an
$n\times n$ matrix. We call $C$ a positive matrix,
and write
$C\geq 0$, if
$c_{ i,j}\geq 0$ for all
$i,j$.  

\medskip A matrix
$A$  is said to be a nonsingular  $M$-matrix   if
\begin{enumerate}
\item[(1)] $a_{ i,j}\leq 0$ for all $i\neq j$.
\item[(2)] $A$ is nonsingular and $A^{ -1}\geq 0$.
\end{enumerate}

It follows from \cite[Lemma 4.2]{EK} that the right-hand side of   (\ref{int.1}) is a Laplace transform for all $\al>0$ if
  $A=K^{-1}$ exists and  is a  nonsingular $M$-matrix.   We refer to   $A$ as the  $M$-matrix  corresponding to $X$. We also use the terminology that $K$ is an inverse $M$-matrix.

   
   Let $K$ be an $n\times n$ matrix with positive entries. We define
\begin{equation}
   K_{Sym}=\{(K_{i,j}K_{j,i})^{1/2}\}_{i,j=1}^{n} .  \label{2.1nn} \end{equation}

 When $A$ is an $M$-matrix, we define 
\begin{equation}
   A_{sym}= \left\{
   \begin{array}{cl}
   A_{j,j}&j=1,\ldots,n\\
  -(A_{i,j}A_{j,i})^{1/2}&i,j=1,\ldots,n, \,\,i\ne j
   \end{array}  \right. .\label{2.15}
   \end{equation}

  In addition, when $A=K^{-1}$ is a non-singular $M$-matrix, it follows from \cite[Lemma 3.3]{MRHD} that $A_{sym}$ is a non-singular  $M$-matrix. We define,
\begin{equation}
     K_{isymi}= (A_{sym})^{-1}.\label{con7}
   \end{equation}
  (The notation $isymi$ stands for, `take the inverse, symmetrize and take the inverse again'.)    
  
  It is obvious that if $K$ is symmetric,
\begin{equation}
     K_{isymi}= K.\label{con7qq}
   \end{equation}
 
  \begin{lemma}\label{lem-2.2} When $K$  is an inverse $M$-matrix so is     $K_{ isymi}$. Furthermore, when,   in addition,   $K$ is     equivalent to a symmetric  inverse  $M$-matrix then,    $K_{Sym}=K_{isymi}$.   
\end{lemma}

\Proof By hypothesis $A=K^{-1}$ is a non-singular $M$-matrix. Therefore, as we just pointed out,  $A_{sym}$ is a non-singular  $M$-matrix. We  denote  it's inverse by $K_{isymi}$. Furthermore, if    $K$ is     equivalent to a symmetric  inverse  $M$-matrix $R$,  then by \cite[Lemma 2.1, $(2.5)$]{MRns} and the fact that $R\geq 0$, we have   $R=K_{Sym}$. 
It is easy to see that the results of \cite[Lemma 2.1]{MRns} hold with $K$ and $R$ replaced by $A $ and $R^{-1}$. Therefore,  since $R^{-1}$ is an $M$-matrix, it follows from the analogue of \cite[Lemma 2.1, $(2.5)$]{MRns} that 
$R^{-1}=A_{sym}$.   Thus $K_{Sym}=K_{isymi}$.   \qed

   The next lemma which is given in \cite[Corollary 3.1]{MRHD}, associates an $\al$-permanental random variable  with kernel $K$ with an $\al$-permanental random variable with a symmetric kernel $K_{isymi}$.  
.


    \begin{lemma}\label{cor-7.1}  For any $\al>0$ let $X_{\al}(n)=(X_{\al,1},\ldots,X_{\al,n})$ be the $\al$-permanental random variable determined by an $n\times n$ kernel $K(n)$ that is an inverse   $M$-matrix and set $A(n)=K(n)^{-1}$. Let $\wt X_\al $   be the $\al$-permanental random  variable determined by     $K(n)_{isymi}$.  Then for all functions $g_{n}$ of $X_\al(n)$ and $\wt X_\al(n)$ and sets $\BB_{n}$ in the range of $g_{n}$
\bea
\lefteqn{ \frac{|A(n)|^{\al}}{|  A(n)_{sym}|^{\al}}  P\(g_{n}( \wt X_\al(n))\in \BB_{n}\)  \le  P\( g_{n}(   X_\al(n))\in\BB_{n}\)\label{3.16w}}\\
 & & \hspace{1in}\le\(1- \frac{|A(n)|^{\al}}{|  A(n)_{sym}|^{\al}} \) +\frac{|A(n)|^{\al}}{|  A(n)_{sym}|^{\al}}P\(g_{n}( \wt X_\al(n))\in\BB_{n}\).\nn
 \eea

 \end{lemma}

  Let $G(1, l,n)$ denote the $n\times n$ matrix obtained by restricting the matrix  $G$ in    (\ref{21.0})  to the $n\times n$ matrix with indices
 $\{1, l+2,\ldots, l+n \}\times \{1, l+ 2 ,\ldots, l+n \}$. 
 Note that $G(1,0,n)=\{G_{i,j};i,j=1,\ldots,n\}:=G(n)$.
 
 Let
    \be K(1,l,n) =\left (
\begin{array}{ ccccc }  1 &1&\ldots&1  \\
1   &G(1,l,n)_{1,1} &\ldots&G(1,l,n)_{1,n}   \\
\vdots& \vdots &\ddots &\vdots  \\
1   &G(1,l,n)_{n,1} &\ldots&G(1,l,n)_{n,n} 
\end{array}\right ). \label{19.19e}
     \ee  
$K(1,l,n)$ is an $(n+1)\times (n+1)$  matrix. We use   $\{0,1\ldots,n\}$ to denote the indices of 
$K(1,l,n)$.

\medskip	We now state three lemmas that are the core of the proof  of  Theorem  \ref{theo-kolexp2q}.   They are proved in Section \ref{sec-isimi}.
 
\bl\label{lem-kmatrix} $K(1,l,n)$ is an inverse M-matrix, and $\AA(1,l,n):=K(1,l,n)^{-1}$ has non-negative row sums. 
\el

\medskip	
The next two lemmas enable  us to use 
Lemma \ref{cor-7.1} in the proof of  Theorem  \ref{theo-kolexp2q}.

\bl\label{lem-koldet} 
\begin{equation}
 1\leq \tau_{1,l,n}:={ | \AA_{sym}(1,l,n) | \over   | \AA(1,l,n) |}\leq 2,\label{22.15}
\end{equation}
and 
\begin{equation}
\lim_{l\to\ff}\tau_{1,l,n}=1,\label{22.15a}
\end{equation}
uniformly in $n$.
\el
 In particular, this shows that we can find some $ \tau_{1,l,\ast}\geq 1$ such that 
 \begin{equation}
 \tau_{1,l,n}\geq \tau_{1,l,\ast} \quad \mbox{ for all $n$ and}\quad\lim_{l\to\ff}\tau_{1,l,\ast}=1.\label{22.15aa}
 \end{equation}

 Set $K_{isymi}(1,l,n) =\( \AA_{sym}(1,l,n)\)^{-1}$.  $\{K_{isymi}(1,l,n)_{j,k};j,k\in  0,1\ldots,n \}$ is an $(n+1)\times (n+1)$  matrix. We use   $\{0,1\ldots,n\}$ to denote the indices of 
$K(1,l,n)_{isymi}$. Let $Z_{(\al)}(n+1)=\{Z_{(\al),0},\ldots, Z_{(\al),n}\}$ be an $\al$-permanental process with kernel $K(1,l,n)$ and $\ov Z_{(\al)}(n+1)=\{\ov Z_{(\al),0},\ldots, \ov Z_{(\al),n}\}$ be an $\al$-permanental process with kernel $K_{isymi}(1,l,n)$. Then it follows from Lemma \ref{cor-7.1}
 that, 
\bea
\lefteqn{ { | \AA_{sym}(1,l,n) | ^{\al}\over   | \AA(1,l,n) |^{\al}}P\(g_{n+1}(  \ov Z_{(\al)}(n+1))\in \BB_{n+1}\)  \le P\(g_{n+1}(    Z_{(\al)}(n+1))\in \BB_{n+1}\) \label{3.16ww}}\nn\\
 & & \le\(1-  { | \AA_{sym}(1,l,n) | ^{\al}\over   | \AA(1,l,n) |^{\al}} \) +{ | \AA_{sym}(1,l,n) | ^{\al}\over   | \AA(1,l,n) |^{\al}} P\(g_{n+1}(  \ov Z_{(\al)}(n+1))\in \BB_{n+1}\) ).\nn
 \eea

  Suppose that $g_{n+1}$ is such that 
\begin{equation}
   g _{n+1}(v_{0}, v_{1},\ldots,v_{n})=g'(v_{1},\ldots,v_{n}).
   \end{equation}
Then 
\begin{equation}
   P\(g_{n+1}(    Z_{(\al)}(n+1))\in \BB_{n+1}\)=P\(g' (    Z_{(\al),1},\ldots, Z_{(\al),n})\in \BB_{n+1}\).
   \end{equation}
and similarly with $Z$ replaced by $\ov Z$. This  is an important observation. We see from (\ref{19.19e}) that $(Z_{(\al),1},\ldots, Z_{(\al),n}))$ is an $\al$-permanental sequence with kernel $G(1,l,n)$. This is what we set out to study. Similarly, $(\ov Z_{(\al),1},\ldots, \ov Z_{(\al),n}))$ is an $\al$-permanental sequence with kernel 
\begin{equation}
     \{K_{isymi}(1,l,n)_{j,k};j,k\in  1\ldots,n \} \label{2.12nn}
   \end{equation}
Since this kernel is symmetric  it is the covariance of a Gaussian sequence
 which we denote by 
$\{\xi(1,l,n)_{j};j\in  1,\ldots,n\}$.

 	 \bl\label{lem-gold}   Let $\{\eta_{j};j\in 1,2,\ldots\} $ be an independent standard normal sequence.   Then,
 \begin{equation}
 \xi(1,l,n)_{1}\stl\sqrt 2 \,\eta_{0}+o_{l}(1) \eta_{0}, \label{22.321q}
 \end{equation}
 and, with $\JJ(n)=\{2,\ldots,n\} $, 
 \begin{equation}
\xi(1,l,n)_{j}\stl \la^{1/2}_{l+j} \,\eta_{ l+j }+\sqrt 2 \,\eta_{0}+o_{l}(1)  \eta_{0}, \hspace{.2 in}j\in\JJ(n),\label{22.32w}
 \end{equation}
and
 \be \{ \xi(1,l,n)_{j}-\xi(1,l,n)_{1};\,\,j\in\JJ(n)\}\label{2.15oo}  \stl \{\la^{1/2}_{ l+j} \,\eta_{ l+j }+o_{l}( \la_{ l+j } )\, \eta_{0};\,\,j\in\JJ(n)\}  .
    \ee 
  \el

 \section{Proof  of  Theorem \ref{theo-kolexp2q}  } \label{sec-3}

 	The results in Theorems  \ref{theo-kolexp2q} are not intuitively obvious. To begin we simplify the problem to show how they come about. We replace $U$ in   (\ref{21.1c})  by  $\ov U=\{\ov U(s,t),s,t\in T\}$ where,  
	\begin{equation}
    \ov U(s,t) =\de_{s,t}\la_{1/s}+2 \mathbf 1 ,\qquad s,t\in T\label{1.56a}, 
    \end{equation} 
  where, $\la_{1/0}:=0$. This is a good approximation to $U$,  since $\lim_{s\to 0}{f_{s}}=\lim_{s\to 0}g_{s}=1$ it is very close to $U$ as $s\to 0$. 
 
  Note that $\ov U$ is the covariance of  a Gaussian sequence $\{\tau_{s},s\in T\}$ where
 \begin{equation}
\tau_{0}=\sqrt 2\eta_{0},\qquad\mbox{and}\qquad   \tau_{s}=\la_{1/s}^{1/2}\eta_{s}+\sqrt 2\eta_{0}, \label{tau}
   \end{equation}
where $\{\eta_{s},s\in T\} $ be an independent standard normal sequence. 
Therefore, $\ov X_{(1/2)}=\{\ov X_{(1/2),s},s\in T\}=\{\tau^{2}_{s}/2,s\in T\}$ is a 1/2 permanental sequence with kernel $\ov U $ and, 
 \begin{equation}
  \{ \ov X_{(1/2),  s}- \ov X_{(1/2),0},s\in T\}= \Big\{\frac{\tau_{s}^{2}  - \tau_{0}^{2}}{2} ,s\in T  \Big\}.\label{3.4mm}
   \end{equation}

 It simplifies things if we use more standard notation. Let $\{\eta_{j};j\in \mathbb N\} $ be an independent standard normal sequence. Define
\begin{equation}
\tau_{0}=\sqrt 2\eta_{0},\qquad\mbox{and}\qquad   \tau_{n}=\la_{n}^{1/2}\eta_{n}+\sqrt 2\eta_{0},\qquad n\ge 2  \label{tau
s}.
   \end{equation}
Then
\begin{equation}
  \{ \ov X_{(1/2),  1/n}- \ov X_{(1/2),0},n\ge 2\}= \Big\{\frac{\tau_{n}^{2}  - \tau_{0}^{2}}{2} ,n\ge 2 \Big\}.\label{3.4mmq}
   \end{equation}
 
 \medskip		Set $\vec \tau _{n}=( \tau _{n,1},\ldots, \tau _{n,k})$ and $\vec \eta _{n}=( \eta _{n,1},\ldots, \eta _{n,k})$ so that,  
    \begin{equation}
     \sum_{i=1}^{k} \tau^{2}_{n,i}=\|  \vec \tau _{n }\|_{2}^{2} \qquad\mbox{and}\qquad     \sum_{i=1}^{k} \eta^{2}_{n,i}=\|  \vec \eta _{n }\|_{2}^{2}\label{3.5mm}.
      \end{equation}
 Consequently, $\ov X_{(k/2)}=\{\ov X_{(k/2),1/n},n\ge 2\}=\{ {\|  \vec \tau _{ n }\|_{2}^{2}}/{2},n\ge 2\} $ is a $k/2$-permanental sequence with kernel $\ov U$,
   and,
\be 
  \{ \ov X_{(k/2),  1/n}- \ov X_{(k/2),0},n\ge 2\}=\Big\{\frac{\|  \vec \tau _{n}\|_{2}^{2}  -  \|  \vec \tau _{0}\|_{2}^{2}}2  ,n\ge 2  \Big\}.\label{2.15mm}
   \ee 
 
 \begin{lemma}\label{lem-3.1} The limits in   (\ref{21.4mva})--(\ref{21.4msa}) hold when $\{X_{(k/2),1/n};n\ge 2\}$ is replaced by $\{\ov X_{(k/2),1/n};n\ge 2\}$.
   \end{lemma}
    
 The proof of Lemma \ref{lem-3.1} uses the next two  lemmas: 

 \begin{lemma} \label{lem-2.6mm} Let   $\eta_{i,j}$, $i=1,\ldots,k$, $j=1,\ldots,n$,  be independent standard normal random variables and let $\vec{\eta  }_{n}=(\eta  _{1,n},\ldots, \eta  _{k,n})$.  Then for any $\vec{v}\in R^{k}$ with   $\|\vec{v}\|_{2}\le 1$,   there exists an ${\Om}'\subset \Om$ with measure 1 such that for all $\om'\in \Om' $ there exists  a subsequence    $\{s_{m}(\om')\}$  of    the integers, such that for all $1\le j\le k$,  
 \begin{equation}
    \lim_{m\to\ff}  \frac{   \eta  _{j, s_{m}(\om')}}{ (2\log s_{m}(\om') )^{1/2}}=  v_{j} ,\qquad a.s. \label{3.20mm}
    \end{equation}
   \end{lemma}

\Proof    Let $\ep_{n}=(\log n)^{-1/k}$. Then,  
\bea
&&   P\(  v_{j}\sqrt 2\log^{1/2} {n}- \ep_{n}  \le \eta  _{j,n}\le   v_{j}\sqrt 2\log^{1/2} {n}+ \ep_{n}  ;j=1,\ldots,k\)\quad \\
&&\qquad=\prod_{j=1}^{k}\nn \frac{1}{\sqrt{2\pi}}\int_{\sqrt 2  v_{j}\log^{1/2} {n}- \ep_{n}  }^{\sqrt 2  v_{j}\log^{1/2} {n}+ \ep_{n}  } e^{-x^{2}/2}\,dx\\&&\qquad\ge\prod_{j=1}^{k}\nn \frac{ \ep_{n}  }{2\sqrt{2\pi}}\displaystyle\frac{1}{n^{   v_{j}^{2}}}\ge{ \ep_{n}^{k}\over (8\pi)^{k/2}  n}=\frac{1}{ (8\pi)^{k/2} n\log n}.
   \eea  
   
The statement in (\ref{3.20mm}) follows from the Borel-Cantelli Lemma.\qed

 	\begin{lemma}\label{lem-3.4mm} Let $\vec{\eta  }_{n}=(\eta  _{1,n},\ldots, \eta  _{k,n})$ be as in Lemma \ref{lem-2.6mm}.  Then, almost surely, the set of limit points of  
\begin{equation}
   \frac{ \vec{\eta}_{n} }{\(\displaystyle 2\log
n \)^{1/2}},  
   \end{equation} 
 is the unit ball in $R^{k}$.   Equivalently,  
 	 \begin{equation}
\sup_{\{\vec w\in R^{k}\,: \, \|\vec w\|_{2}\leq 1\}}\liminf_{n\to\ff}\Big\|\frac{ \vec{\eta}_{n} }{\(\displaystyle 2\log
n \)^{1/2}} - \vec w\,\Big\|_{2}=0,\qquad a.s.\label{Note.30}
\end{equation}
\end{lemma}

\Proof   It follows from Lemma \ref{lem-2.6mm} that 
	 \begin{equation}
\liminf_{n\to\ff}\Big\|\frac{ \vec{\eta}_{n} }{\(\displaystyle 2\log
n \)^{1/2}} - \vec w\,\Big\|_{2}=0,\qquad a.s.\label{Note.30j}
\end{equation}
 holds for any fixed $\vec w\in R^{k}$ with   $\|\vec w\|_{2} \le 1$. Furthermore, since (\ref{Note.30j}) holds almost surely, it can be extended to hold on a countable dense subset of the unit ball of $R^{k}$. We can then extend it to all $\vec w $ in the unit ball of $R^{k}$.    \qed

 	\noindent {\bf  Proof of  Lemma \ref{lem-3.1}} 
We write,
    \bea
    \frac{ \|\vec \tau_{n}  \|_{2}^{2}  - \|\vec\tau_{0}  \|_{2}^{2}}{2}  &=&
  \frac{   \|\vec \tau_{n} -\vec \tau_{0} \|_{2}^{2}}{2}+ (\vec \tau_{n} -\vec \tau_{0} )\cdot \vec \tau_{0} \nn\\
     &=&
  \frac{   \la_{n} \|\vec\eta_{n}  \|_{2}^{2}}{2}+ \sqrt{2\la_{n}}(\vec\eta_{n}\cdot\vec\eta_{0} ) .\label{2.23mmo}
     \eea	Consequently,
 \be 
  \frac{   \|\vec \tau_{n}  \|_{2}^{2}  - \|\vec\tau_{0}  \|_{2}^{2} }{2(\la_{n}\log n)^{1/2} } =
   \frac{ (\la_{n}\log n)^{1/2}  \|\vec\eta_{n}  \|_{2}^{2}}{2\log n}+  \frac{ {2  }(\vec\eta_{n}\cdot\vec\eta_{0} ) }{(2\log n)^{1/2}}    .\label{3.23kk}
     \ee  
Let $\vec u$ be a random variable in the unit ball of $R^{k}$ and consider,
\be \bb^{1/2} \|\vec u\|_{2}^{2}+2    (\vec u\cdot\vec\eta_{0} ).  \label{3.24oo}
   \ee 
It is clear that to maximize this term we should take $\vec u$ in the direction of  $\vec \eta_{0}$ with norm 1. That is we should take $\vec u=\vec \eta_{0}/\|\vec \eta_{0}\|_{2}$, in which case (\ref{3.24oo}) is equal to,
\be \bb^{1/2}  +2   \|\vec \eta_{0}\|_{2}.  \label{3.24jj}
   \ee 
It now follows from Lemma \ref{lem-3.4mm} with $\vec w=\vec \eta_{0}/\|\vec \eta_{0}\|_{2}$ that, 
\begin{equation}
   \limsup_{n\to\ff} \frac{   \|\vec \tau_{n}  \|_{2}^{2}  - \|\vec\tau_{0}  \|_{2}^{2} }{2(\la_{n}\log n)^{1/2} } =\bb^{1/2}  +2   \|\vec \eta_{0}\|_{2}.    \end{equation}
    Therefore, it follows from   (\ref{2.15mm}) that (\ref{21.4m}) holds when $\{X_{(k/2),1/n} \}$ is replaced by $\{\ov X_{(k/2),1/n} \}$. A trivial modification of this argument shows that (\ref{21.4mv}) also  holds  when  $\{X_{(k/2),1/n} \}$ is replaced by $\{\ov X_{(k/2),1/n} \}$.

\medskip	Consider (\ref{3.24oo}) again and assume that $\bb>0$.     It is clear that   for any fixed  value of $\|\vec u\|_{2} $, (\ref{3.24oo}) is minimized when  ${\vec u} $   in the     negative of the direction of  ${\eta_{0}} $. Let ${\vec u}  =-\ga {\vec \eta_{0}} $ where $0<\ga\le 1/\| \vec\eta_{0}\|_{2} $; (so that $\| \vec u \|_{2}\le 1 $). With this substitution   (\ref{3.24oo}) is,
 \begin{equation}
  \(  \bb^{1/2}\ga^{2} -2  \ga \) \| \vec\eta_{0}\|_{2}^{2}.  \label{3.25mm}
    \end{equation}
 This is minimized when $\ga=1/\bb^{1/2}$. For this value of $\ga$, ${\vec u}  =-  {\vec \eta_{0}}/\bb^{1/2} $. Therefore,  as long as $\| \vec\eta_{0}\|_{2}\le \bb^{1/2}$ we can achieve the minimum which is $-\| \vec\eta_{0}\|_{2}^{2}/\bb^{1/2}$.
  
   When  $\| \vec\eta_{0}\|_{2}>\bb^{1/2}$ we are restricted to $0<\ga\le 1/\bb^{1/2}$. One can check that $\bb^{1/2}\ga^{2} -2  \ga$ is decreasing for $\ga$ in this range. Therefore it takes its minimum at $\ga=1/\| \vec\eta_{0}\|_{2}$. For this value of $\ga$, (\ref{3.25mm}) is equal to $\bb^{1/2}-2 \| \vec\eta_{0}\|_{2}$.      Consequently,
 \be  \liminf_{s\to 0}\frac{   \|\vec \tau_{s}  \|_{2}^{2}  - \|\vec\tau_{0}  \|_{2}^{2} }{2 (\la_{s}\log 1/s)^{1/2}}\label{21.4pp} =  \left\{\begin{array} {cc}  -   \| \vec\eta_{0}\|_{2}^{2}/\bb^{1/2}, & \mbox{if }\hspace{.2 in}   \| \vec\eta_{0}\|_{2}^{2}< \bb\\
\bb^{1/2} -2    \| \vec\eta_{0}\|_{2} ,& \mbox{if }\hspace{.2 in}  \| \vec\eta_{0}\|_{2}^{2}>\bb
\end{array}\right.  
\ee

 When $\bb=0$ the expression in  (\ref{3.24oo}) is $2(\vec u\cdot\vec\eta_{0})$. 
 This is minimized when $\vec u=-\eta_{0}/ \| \vec\eta_{0}\|_{2}$, keeping in mind that we must have $\vec u$ in the unit ball of $R^{k}$. With this substitution
 $2(\vec u\cdot\vec\eta_{0})=-2 \| \vec\eta_{0}\|_{2}$. Therefore, (\ref{21.4pp}) also holds when $\bb=0.$
Using  (\ref{2.15mm}) we see that (\ref{21.4ms}) holds when when $\{X_{(k/2),1/n} \}$ is replaced by $\{\ov X_{(k/2),1/n} \}$.\qed

We now establish the material we need to use Lemma \ref{cor-7.1}. 
On page  \pageref{2.12nn}  we define the Gaussian sequence $\{\xi(1,l,n)_{j};j\in  1,\ldots,n\}$ with covariance $K_{isymi}(1,l,n)=  \{K_{isymi}(1,l,n)_{j,k};j,k\in  1\ldots,n \}$. Let   $\{\xi(1,l,n)_{j,i};j\in  1,\ldots,n; i=1,\ldots,k\}$ be $k$ independent copies of $\{\xi(1,l,n)_{j};j\in  1,\ldots,n\}$. Then 
\begin{equation}
Z(1,l,n)_{(k/2),j}  =  \sum_{i=1}^{k} { \xi^{2}(1,l,n)_{j,i}\over 2},\qquad j=1,\ldots,n,\label{3.21}
   \end{equation}
is a $k/2$-permanental sequence with kernel $K_{isymi}(1,l,n)$.

\begin{lemma} Let $Z(1,l,n)_{(k/2) } =\{Z(1,l,n)_{(k/2),j};j=1,\ldots,n \}$ and define $\ov X(1,l,n)_{(k/2) }$ similarly. Let $\JJ(n)=\{2,\ldots,n\} $. Then,  
 \bea
 \lefteqn{\bigg\{\frac{Z(1,l,n)_{(k/2),j }-Z(1,l,n)_{(k/2),1}}{2}\,;j\in\JJ(n)\bigg\}\label{3.28mm}}\\
&& =  
\bigg\{\frac{ { \ov X_{(k/2),  1/(l+j)}- \ov X_{(k/2),0}}}{2}   +  o (  \la  _{ l+j}  ) \( \|  \ov X_{(k/2),0}\|^{2}_{2}+1\);\,j\in\JJ(n)\bigg\}\nn.
     \eea 

 \end{lemma}
 
 \Proof
To begin we    show  that,  
 \bea
 &&\Big\{  {    \xi^{2} (1,l,n)_{j}- \xi ^{2}(1,l,n)_{1}  };j\in\JJ(n) \Big\}\label{3.27mm}\\
 &&\qquad =  \Big\{
 { \tau_{ l+j  }^{2}  - \tau^{2}_{0}  } + o( \la_{ l+j} )\,      ( \tau ^{2}_{0}+1) ;j\in\JJ(n)  \Big\}.\nn
     \eea 
   It follows from    (\ref{2.15oo}) that,     
      \bea
 \lefteqn{ \hspace{-.075in}  \Big\{ (\xi (1,l,n)_{j}- \xi (1,l,n)_{1})^{2} ;j\in\JJ(n) \Big\}  } \\
     && \hspace{-.2in} = \Big\{\la _{ l+ j} \,\eta^{2}_{ l+ j}+\la^{1/2}_{ l+ j} \,\eta_{ l+ j}o( \la_{ l+ j} )\, \eta_{0}\nn+o( \la^{2}_{ l+ j} )\, \eta^{2}_{0};j\in\JJ(n) \Big\}.
        \eea
 Using this (\ref{22.321q}) and (\ref{2.15oo}) and the relationship $a^{2}-b^{2}=(a-b)^{2} +2(a-b)b$,  we see that,  
   \bea      {    \xi^{2} (1,l,n)_{j}- \xi ^{2}(1,l,n)_{1} } &=&  { \la _{ l+ j}\eta^{2}_{ l+ j}} +  { 2 \sqrt2 \la^{1/2} _{ l+ j}\eta _{ l+ j}}\eta_{0}  +\psi(1,l,n)_{j} \nn\\
      &=&  { \tau_{ l+ j }^{2}  - \tau^{2}_{0}  } +\psi(1,l,n)_{j} \nn,
      \eea        
   where
   \begin{equation}
      \psi(1,l,n)_{j}=    \la^{1/2}_{ l+ j} \,\eta_{ l+ j}o( \la_{ l+ j} )\, \eta_{0}\nn+o( \la^{2}_{ l+ j} )\, \eta^{2}_{0}
 + {o( \la  _{ l+ j} )\, \eta^{2} _{0} } .
      \end{equation}  
  Since,   
   \begin{equation}
        \la^{1/2}_{ l+ j} \,\eta_{ l+ j}\eta_{0}=\frac{     (\la_{ l+ j}\log(l+ j))^{1/2} \,\eta_{ l+ j}}{(\log   (l+ j))^{1/2}}\eta_{0},
      \end{equation} 
        it follows from  (\ref{21.3b})     that,  
 \begin{equation}
    \limsup_{l\to\ff} \la^{1/2}_{ l+ j} \, \eta_{ l+ j} \eta_{0}  =(2\bb)^{1/2} \limsup_{l\to\ff}\frac{     \,\eta_{ l+ j}}{(2\log   (l+ j))^{1/2}}\eta_{0}=(2\bb)^{1/2}|\eta_{0}|\label{5.11oyww},\qquad a.s. 
    \end{equation}    
    where we use the well known facts that 
    \begin{equation}
       \limsup_{l\to\ff}\frac{   \pm   \,\eta_{ l+ j}}{(2\log   (l+ j))^{1/2}}=1,\qquad a.s.
       \end{equation}   
       and is independent of $\eta_{0}$.        
   Consequently,   
   \begin{equation}
      \psi(1,l,n)_{j}=   o( \la  _{ l+ j} )\,  ( \eta ^{2}_{0} +|\eta_{0}|).    \label{3.28}  \end{equation} 
This gives (\ref{3.27mm}).  
   
    Using (\ref{3.27mm})
and  
 (\ref{2.15mm}) we get,  
  \bea
  && \sum_{i=1}^{k} { \xi^{2}(1,l,n)_{j,i}- \xi^{2}(1,l,n)_{1,i}\over 2}\label{3.28mmr}\\
  &&\qquad=  
\frac{\|  \vec \tau _{ l+j}\|_{2}^{2}  -  \|  \vec \tau _{0}\|_{2}^{2}}2 +   o( \la  _{l+j} )(\|   \vec \tau _{0}\|_{2}^{2}+1)\nn\\
&&\qquad=  
 { \ov X_{(k/2), 1/( l+j)}- \ov X_{(k/2),0}}  +  o( \la  _{l+j} )  \(\|  \ov X_{(k/2),0}\|^{2}_{2}+1\)\nn,
     \eea 
or equivalently (\ref{3.28mm}). \qed

 	\noindent{\bf   Proof  of  Theorem  \ref{theo-kolexp2q}   } 	  
 Let 
 \begin{equation}
    \Phi=\Phi(\bb,X_{(k/2),0})= \bb^{1/2}   +2X_{(k/2),0}  ^{1/2}
    \end{equation}It follows from Lemmas \ref{cor-7.1} applied to the left-hand side of (\ref{3.28mm}) and  Lemma \ref{lem-koldet}  that,
\begin{eqnarray}
&&P\(\Big|\sup_{ l+2\leq  j\leq l+n}  \frac{   X_{(k/2),  1/j}-   X_{(k/2),0}}{(\la_{j}\log j)^{1/2}} -\Phi\Big |\leq \ep\)
\label{first} \\
&&\nn\qquad \geq    \tau^{-k/2}_{1,l,n} P\( \Big|\sup_{ l+2\leq  j\leq l+n}  \frac{\ov   X_{(k/2),  1/j}- \ov   X_{(k/2),0}}{(\la_{j}\log j)^{1/2}}- \Phi\right. \\
&&\hspace{2 in}\left.  + o\( \la_{j}^{1/2}/(\log  j)^{1/2} \)\(  \|  \ov X_{(k/2),0}\|^{2}_{2}+1)\)\Big |\leq \ep \)
 \nn\\
&&\nn\qquad \geq    \tau^{-k/2}_{1,l,n} P\(  \Big|\sup_{ l+2\leq  j\leq l+n}  \frac{\ov   X_{(k/2),  1/j}- \ov   X_{(k/2),0}}{(\la_{j}\log j)^{1/2}}- \Phi\Big |\right. \\
&&\hspace{3.1in}  \leq \ep/2 \Bigg)- o\(  1/(\log l) ^{1/2}\) ,
 \nn
\end{eqnarray}
where we use  the inequality $P(A\cap B)$=$P(A)-P(A\cap B^{c})\ge P(A)-P(  B^{c})$ to get  the last line.

 	 We take the limit as  
 $n\to\ff$, and use (\ref{22.15aa}), to obtain,  
 \bea &&  P\(\Big|\sup_{j\ge  l+2 }  \frac{   X_{(k/2),  1/j}-   X_{(k/2),0}}{(\la_{j}\log j)^{1/2}} -\Phi\Big |\leq \ep\)
\label{22.45mm} \\
&&\nn\qquad \geq    \tau^{-k/2}_{1,l,*}  P\(  \Big|\sup_{ j\ge l+2 }  \frac{\ov   X_{(k/2),  1/j}- \ov   X_{(k/2),0}}{(\la_{j}\log j)^{1/2}}- \Phi\Big |\right. \\
&&\hspace{3.1in}  \leq \ep/2 \Bigg)- o\(  1/(\log l) ^{1/2}\)  .
\label{22.45}\nn
\end{eqnarray}  
  Using (\ref{22.15aa}) again and  taking the limit as  
 $l\to\ff$ we have, 
 \bea && P\(\Big|\limsup_{   j\to\ff}  \frac{   X_{(k/2),  1/j}-   X_{(k/2),0}}{(\la_{j}\log j)^{1/2}} -\Phi\Big |\leq \ep\)
\label{22.45mma} \\
&&\nn\qquad \geq    P\(\Big| \limsup_{   j\to\ff}  \frac{\ov   X_{(k/2),  1/j}-   \ov X_{(k/2),0}}{(\la_{j}\log j)^{1/2}} -\Phi\Big |   \leq \ep/2 \)  =1 .
\label{22.45x}\nn
\end{eqnarray}  
where for the last inequality we use Lemma \ref{lem-3.1}. This gives(\ref{21.4m}).

The proofs of (\ref{21.4mv}) and  (\ref{21.4ms})  proceed in exactly the same way, beginning with putting the appropriate term in the first line of (\ref{first}).\qed

\section{  Proof that  Theorem \ref{theo-kolexp2q} and  Theorem \ref{theo-1.4}  are equivalent}\label{sec-4}
We first show that Theorem \ref{theo-kolexp2q} implies Theorem \ref{theo-1.4}. We point out following the statement of  Theorem \ref{theo-kolexp2q} that it is a restatement of Theorems \ref{theo-kolexp12} and (\ref{theo-kolexp22}). We use results in Theorems \ref{theo-kolexp12} and (\ref{theo-kolexp22}) to obtain Theorem \ref{theo-1.4}. 

When  $\bb>0$ we can write (\ref{21.4m}) as,
\be
 \limsup_{ n\to\ff  } { X _{(k/2), 1/n} }  =  ( X^{1/2}_{(k/2),0}+\bb^{1/2})^{2} \,, \qquad a.s.\label{zz1}
\ee
Taking square roots we get, 
\begin{equation}
    \limsup_{  n\to\ff } { X^{1/2} _{(k/2), 1/n} }    =  X^{1/2}_{(k/2),0}+\bb^{1/2}   \,, \qquad a.s.\label{zz2}
   \end{equation}
 This gives both (\ref{firstqq})  and    (\ref{firstq}).

 To consider the case in which $\bb=0$  we 
  write,
\be 
  \frac{ X  _{(k/2), 1/n}- X _{(k/2),0} }{(\la_{n}\log n)^{1/2}} =\frac{ \(X^{1/2}_{(k/2), 1/n}- X^{1/2}_{(k/2),0}\) \(X^{1/2}_{(k/2), 1/n}+X^{1/2}_{(k/2),0}\)}{(\la_{n}\log n)^{1/2}} \label{1.24mm}
 \ee
and note that when $\bb=0$ it follows from (\ref{21.4n}) that $X^{1/2}_{(k/2), 1/n}$ is continuous at 0. Consequently, 
\begin{equation}
   \limsup_{n\to\ff} \frac{ X  _{(k/2), 1/n}- X _{(k/2),0} }
   {(\la_{n}\log n)^{1/2}} = \limsup_{n\to\ff} \frac{ X^{1/2}  _{(k/2), 1/n}- X^{1/2} _{(k/2),0} }
   {(\la_{n}\log n)^{1/2}}2X _{(k/2),0}^{1/2}   \end{equation}
Therefore,   (\ref{21.4nq}) implies (\ref{firstq}). A similar argument gives (\ref{firstqq}).

\medskip	When $X_{(k/2),0}\ge \bb$ we see from (\ref{21.4ms}) that
\begin{equation}
  \liminf_{ n\to\ff} { X_{(k/2), 1/n}  } =(X^{1/2}_{(k/2),0} -\bb^{1/2} )^{2}. 
\end{equation}
  If $\bb>0$ we take the square root and rearrange this to get (\ref{secondq}).

  If $\bb=0$ we see by (\ref{21.4nn}) and the fact that $\lim_{n\to\ff}X_{(k/2), 1/n}=X_{(k/2), 0}$ in this case that,
\begin{equation}
\liminf_{  n\to \ff }{X^{1/2}_{(k/2), 1/n}- X^{1/2}_{(k/2),0} \over  (\la_{n}\log n)^{1/2}}2 X^{1/2}_{(k/2),0} =-2 X^{1/2}_{(k/2),0}\,\,, \hspace{.1 in}a.s.\label{dd}
\end{equation}
This gives (\ref{secondq}) when $\bb=0$.
 
 When $X_{(k/2),0}< \bb$, which implies, in particular, that $\bb>0$, it follows from (\ref{21.4ms}) that
   \begin{equation}
     \liminf_{n\to\ff }  X^{1/2}_{(k/2), 1/n}   =0.
     \end{equation}
 This gives  (\ref{secondq}) when $X_{(k/2),0}< \bb$.

\medskip	We now show that Theorem \ref{theo-1.4} implies Theorem \ref{theo-kolexp2q}. Suppose   (\ref{firstq}) holds and $\bb>0$. Then 
\begin{equation}
\limsup_{n\to\ff}   X^{1/2}_{(k/2), 1/n}=X^{1/2}_{(k/2),0}+\bb^{1/2}.
   \end{equation}
or equivalently,  
\begin{equation}
  \limsup_{n\to\ff}  X _{(k/2), 1/n}-X_{(k/2),0}=\bb +2\bb^{1/2}X^{1/2}_{(k/2),0} .
   \end{equation}
This gives (\ref{21.4m}) when $\bb>0$.

Suppose   (\ref{firstq}) holds and $\bb=0$. This implies that $X  _{(k/2), 1/n}$ is continuous at 0. We write,
\be 
  \frac{ X  _{(k/2), 1/n}- X _{(k/2),0} }{(\la_{n}\log
n)^{1/2}} =\frac{ \(X^{1/2}_{(k/2), 1/n}- X^{1/2}_{(k/2),0}\) \(X^{1/2}_{(k/2), 1/n}+X^{1/2}_{(k/2),0}\)}{(\la_{n}\log
n)^{1/2}} \label{1.24mmx}
 \ee
which by  (\ref{firstq}), gives
\begin{equation}
   \limsup_{n\to\ff} \frac{ X  _{(k/2), 1/n}- X _{(k/2),0} }{(\la_{n}\log
n)^{1/2}} =2X _{(k/2),0}^{1/2}   \end{equation}
This gives (\ref{21.4nq}) when $\bb=0$.  

\medskip	Suppose   (\ref{secondq}) holds and $ {X_{(k /2),0}^{1/2}}/{\bb}<1 $, which implies, in particular that $\bb>0$. Therefore,
\begin{equation}
     \liminf_{ n\to\ff }\frac{ X^{1/2}_{(k/2), 1/n}  }{(\la_{n}\log
n)^{1/2}}=0,
   \end{equation}
or, equivalently
\begin{equation}
     \liminf_{n\to\ff } { X _{(k/2), 1/n}  } =0,
   \end{equation}
since $ \bb >0$. This gives  (\ref{21.4ms}) when $ {X_{(k /2),0}^{1/2}}/{\bb}<1 $.

Now suppose   (\ref{secondq}) holds, $ {X_{(k /2),0}^{1/2}}/{\bb}\ge 1 $, and $\bb>0$. We have,
\be     \liminf _{n\to\ff }\frac{ X^{1/2}_{(k/2), 1/n}  }{(\la_{n}\log
n)^{1/2}} = \frac{X^{1/2}_{(k/2),0}}{\bb^{1/2}   }-1,\ee 
or, equivalently,
\be    \liminf _{n\to\ff } { X^{1/2}_{(k/2), 1/n}  }   =   X^{1/2}_{(k/2),0} - \bb^{1/2}  .\ee 
Squaring both sides and rearranging, we get
  (\ref{21.4ms}) when $ {X_{(k /2),0}^{1/2}}/{\bb}\ge 1 $.

When $\bb=0$, $X _{(k/2), 1/n}$ is continuous at 0. Therefore,
\bea  &&\liminf_{ n\to\ff}\frac{ X _{(k/2), 1/n}- X _{(k/2),0} }{(\la_{n}\log
n)^{1/2}}\\&&\qquad=\liminf_{ n\to\ff}\frac{ (X^{1/2}_{(k/2), 1/n}- X^{1/2}_{(k/2),0}) (X^{1/2}_{(k/2), 1/n}+X^{1/2}_{(k/2),0})}{(\la_{n}\log
 n)^{1/2}}\nn\label{1.23mqm}\\
&&\qquad=\liminf_{ n\to\ff}\frac{ (X^{1/2}_{(k/2), 1/n}- X^{1/2}_{(k/2),0}) 2X^{1/2}_{(k/2),0} }{(\la_{ n}\log
 n)^{1/2}}.\nn
\eea
Using (\ref{secondq}) this gives (\ref{21.4nn})  when $\bb=0$. \qed

\section{   Proofs of  Theorems   \ref{theo-1.1mm},  \ref{theo-1.2mm} and \ref{theo-kolexp2s} }   \label{sec-5}

\noindent {\bf  Proof  of  Theorem  \ref{theo-1.1mm} }	
Lemma \ref{lem-kmatrix} 
 holds for $l=0$. As explained in the beginning of Section \ref{sec-2}, 
it then follows    that for any $n$,   $G(n)$ is the kernel of an $\al$-permanental vector for each $\al>0$. Therefore,   by the Kolmogorov extension theorem we have  that for any $\al>0$ there exists a permanental process $X_{(\al)}=(X_{(\al),1},X_{(\al),2}, \ldots)$ with kernel $G$. (In Section \ref{sec-1} we write this as   $X_{(\al)}=(X_{(\al),0},X_{(\al),1/2} ,\ldots)$.)\qed

\noindent {\bf  Proof  of  Theorem  \ref{theo-kolexp2s} }   In the notation of \cite[Theorem 1.6]{MRsp}
  \begin{equation}
\si^{2}(1/n,0)=\la_{n}+1+f _{n}g_{n}+2-2\((1+f_{n})(1+g_{n})\)^{1/2}.
     \end{equation}
  Suppose that $g_{n}\le f_{n}$. Then
  \begin{equation}
  \la_{n}+4
\ge  \la_{n}+1+f _{n}g_{n}+2\ge \la_{n}+4-(1-g_{n}^{2})  .      \end{equation}
  and
  \begin{equation}
4\ge    2\((1+f_{n})(1+g_{n})\)^{1/2}\ge 2(1+g_{n})  =4-2(1-g_{n}) .
     \end{equation}
  Therefore,  
  \begin{equation}
     \si^{2}(1/n,0)\ge  \la_{n}+4-(1-g_{n}^{2}) -4=  \la_{n}-2(1-g_{n}),
     \end{equation}
  and
  \begin{equation}
      \si^{2}(1/n,0)\le  \la_{n} +2(1-g_{n}).
     \end{equation}
  It follows from (\ref{21.2}) that $(1-g_{n})=o(\la_{n})$. Consequently,
    \begin{equation}
      \si^{2}(1/n,0)=  \la_{n}(1+o(1)).\label{1.30}
     \end{equation}
     The same argument gives (\ref{1.30}) when $f_{n}\le g_{n}$.
  Therefore,   $\la_{n}\log n \to \bb$ implies that $  \si^{2}(1/n,0)\log n \to \bb$. 
   The theorem now follows from \cite[Theorem 1.6]{MRsp}. \qed
  
  \vspace{-.05in}
   \noindent{\bf Proof of Theorem \ref{theo-1.2mm} } \label{page-theo-1.2}  It follows from (\ref{1.30}) that  $\la_{n}\log n \to \ff$ implies that $  \si^{2}(1/n,0)\log n \to \ff$. 
  Consequently the upper bound in (\ref{21.4}) follows from \cite[Theorem 1.6]{MRsp}. The lower bound follows from (\ref{22.6}) and \cite[Lemma 7.3]{MRsp}.\qed
  
\vspace{-.1 in} \begin{remark} {\rm 
 We can prove almost all of Theorem \ref{theo-1.2mm} by the methods of this paper. By (\ref{3.23kk}),
\be 
  \frac{   \|\vec \tau_n  \|_{2}^{2}  - \|\vec\tau_{0}  \|_{2}^{2} }{2 \la_n\log n  } =
   \frac{    \|\vec\eta_n  \|_{2}^{2}}{2\log n}+  \frac{ {2  }(\vec\eta_n\cdot\vec\eta_{0} ) }{ 2\log n }    .\label{ 4.7}
     \ee  
Therefore, by (\ref{3.5mm}) and (\ref{2.15mm}), for all $k\ge 1$
 \be
\limsup_{n\to\ff}\frac{ \ov X_{(k/2), 1/n}-\ov X_{(k/2), 0}}{  \la_n\log
n  }=1\qquad 
a.s.\label{21.4oo}
\ee  
 Since  $\lim_{n\to\ff }  \la_n\log n  =\ff$, this is the same as,
\be
\limsup_{n\to\ff}\frac{ \ov X_{(k/2), 1/n} }{  \la_n\log
n  }=1\qquad 
a.s.,\qquad \forall k\ge 1\label{21.4pps}
\ee
We now use the facts that for all $\al>0$, $\{   X_{(\al), 1/n};n\ge 2\}$  is positive and infinitely divisible to see that,
\be
\limsup_{n\to\ff}\frac{ \ov X_{(\al), 1/n} }{  \la_n\log
n  }\le 1\qquad 
a.s.,\qquad \forall \al>0,\label{21.4qq} 
\ee
and
\be
\limsup_{n\to\ff}\frac{ \ov X_{(\al), 1/n} }{  \la_n\log
n  }\ge 1\qquad 
a.s.,\qquad \forall \al>1/2.\label{21.4qqx} 
\ee
As in the proof of Theorem \ref{theo-kolexp2q} it is easy to extend these results to $\{  X_{(\al), 1/n};n\ge 2\}$.\qed
 
  }\end{remark}

\section{Proofs of Lemmas \ref{lem-kmatrix},    \ref{lem-koldet} and \ref{lem-gold}}\label{sec-isimi}
The proofs are rather formal. We begin by trying to explain why we make the substitutions that lead to them.  Our goal is to  use Lemma \ref{cor-7.1} to get probability estimates for the $k/2$-permanental sequences $ (X_{(k/2) ,1} ,X_{(k/2) ,l+2},   \ldots, \newline	  X_{(k/2) ,l+n})$ with kernel $G(1,l,n) $ defined just above  (\ref{19.19e}). To use Lemma \ref{cor-7.1} we must be able to find $G^{-1}(1,l,n) $. This is not easy to do. We accomplish this by embedding  $G(1,l,n) $ in $K(1,l,n) $ as in (\ref{19.19e}). We use Lemma \ref{lem-7.1}  to find the inverse of $K(1,l,n) $.   

 	\medskip	\noindent{\bf  Proof of Lemma \ref{lem-kmatrix} }  
Let 
\be
H(1,l,n) =G(1,l,n)-{\bf  1},\label{sec5.0}
\ee
  where ${\bf  1}$ is the $n\times n$ matrix with all entries equal to $1$.
Therefore,
        \be H(1,l,n)  =\left (
\begin{array}{ cccc  } 1 & g_{l+2 }  &\dots& g_{l+n } \\
  f_{l+2 }&\la_{l+2 }+f_{l+2 } g_{l+2 }& \dots &f_{l+2 }g_{l+n }  \\
\vdots&\vdots& \ddots&\vdots \\
 f_{l+n } &f_{l+n }g_{l+2 }&  \dots &\la_{l+n }+f_{l+n }g_{l+n }   \end{array}\right ).\label{22.5}
     \ee 
     By subtracting $f_{l+j}$ times the first row from the $j$-th row for $j=2,\ldots, n$, we see that,
      \begin{equation}
| H(1,l,n) | =\prod_{j =l+2 }^{l+n }\la_{j}.   \label{22.10}
 \end{equation}
   Note that,  
      \be      H(1,l,n)^{-1}=\left (
\begin{array}{ cccc  }  (1+\sum_{j =l+2 }^{l+n } f_{j}g_{j}/\la_{j})&-g_{l+2}/\la_{2+l} &\dots&- g_{l+n}/\la_{l+n}  \\
-f_{l+2}/ \la_{l+2}&1/\la_{l+2}   & \dots &0   \\
\vdots&\vdots& \ddots&\vdots \\
-f_{l+n}/ \la_{l+n} &0 &  \dots &1/\la_{l+n}     \end{array}\right ).\label{22.6}
    \ee
    This is easy to verify by computing $H(1,l,n)H(1,l,n)^{-1}$. Nevertheless, it is useful to see how we obtain $H(1,l,n)^{-1}$. To simplify the notation a little 
 we do this when $l=0$ and set $H(1,0,n)=H(n)$. We note that
 \begin{equation}
    H(n)=D_{\bf f}JD_{\bf g}
    \end{equation}
 where ${\bf f}=(1,f_{2},\ldots,f_{n})$, ${\bf g}$ is defined similarly and,
   \begin{equation}
      J=D_{\mathbf  r}+\bf 1=\left (
\begin{array}{ ccccc }1  &{\bf  1} \\
{\bf  1}^{T}  &D_{\mathbf  r'}  \end{array}\right ).
      \end{equation}
where ${\bf r}=(1,\la_{2}/f_{2}g_{2},\ldots,\la_{n}/f_{n}g_{n})=:(1,r_{2},\ldots,r_{n})$ and ${\bf r'} =(r_{2},\ldots,r_{n})$. 

It follows from Lemma \ref{lem-7.1} that
 \be J^{-1}  =\left (
\begin{array}{ ccccc }  1+\sum_{j =1}^{n} \ds\frac{f_{j}g_{j}}{\la_{j}}  &-\ds\frac{f_{1}g_{2}}{\la_{2}} &\ldots&-\ds\frac{f_{1}g_{n}}{\la_{n}}   \\&&\vspace{-.1 in}\\
-\ds\frac{f_{2}g_{1}}{\la_{2}}   &\ds\frac{f_{2}g_{2}}{\la_{2}} &\ldots&0 \\
\vdots& \vdots &\ddots &\vdots  \\
-\ds\frac{f_{n}g_{1}}{\la_{n}}    &0 &\ldots&\ds\frac{f_{n}g_{n}}{\la_{n}}  
\end{array}\right ). \label{19.ll}
     \ee  
 We get (\ref{22.6}) by taking
 \begin{equation}
    H^{-1}(n)=D_{\bf 1/f}J^{-1}D_{\bf 1/g}
    \end{equation}
where  ${\bf 1/f}=(1,1/f_{2},\ldots,1/f_{n})$ and ${\bf 1/g}$ is defined similarly.

\medskip	We now find the inverse of $K(1,l,n)$. Consider $K(1,l,n)$ in (\ref{19.19e}). One can see from   (\ref{sec5.0}) that, 
\be
G(1,l,n)= H(1,l,n)+\bf 1 .\label{5.9mm}
\ee
Therefore, it follows from Lemma \ref{lem-7.1} with $h_{k}=1$, $k=1,\ldots,n$, that
   
  \be  K(1,l,n)^{-1}  =\left (
\begin{array}{ cccc  } 1+\rho(1,l,n)&- c(1,l,n)_{1 }  &\dots&- c(1,l,n)_{n } \\
 - r(1,l,n)_{1 }&H(1,l,n)^{1,1} & \dots &H(1,l,n)^{ 1,n}   \\
\vdots&\vdots& \ddots&\vdots \\
  -r(1,l,n)_{n } &H(1,l,n)^{n,1}&  \dots &H (1,l,n)^{n,n}     \end{array}\right ),\label{22.18}
     \ee 
     where
  \bea
r(1,l,n)_{1 }=   \sum_{k=1}^{n} H(1,l,n)^{1 ,k }&=&1-   \sum_{j =l+2 }^{l+n } g_{j}\(1-f_{j}\)\frac{1}{\la_{j}}, \label{22.11j}
\\
 r(1,l,n)_{j }=    \sum_{k=1}^{n} H(1,l,n)^{j ,k }&=&\(1- f_{l+j} \)\frac{1}{\la_{l+j}},\hspace{.2 in}j=2,\ldots,n,\label{22.12j}\\
   c(1,l,n)_{1 }=   \sum_{j=1}^{n} H(1,l,n)^{j ,1 }&=&1-   \sum_{j =l+2 }^{l+n }f_{j} \( 1-g_{j}\)\frac{1}{\la_{j}},\label{22.13j}\\
 c(1,l,n)_{k }=   \sum_{j=1}^{n} H(1,l,n)^{j ,k }&=&\(1-g_{l+k} \)\frac{1}{\la_{l+k}},\hspace{.2 in}k=2,\ldots,n.\label{22.14j}
 \eea
and
  \begin{equation}
\rho(1,l,n)= \sum_{j,k=1}^{n}H(1,n,l)^{j ,k }=1+  \sum_{j =l+2 }^{l+n } \(1- g_{j} \)\(1-f_{j}\)\frac{1}{\la_{j}}.\label{22.19}
 \end{equation}
 
 Note that $\rho1,l,n)$ is the sum of all the elements in $H(1,l,n)^{-1}$, 
$r(1,l,n)_{j }$ is the sum of the  $j$-th row of $H(1,l,n)^{-1}$       and  
$c(1,l,n)_{k }$ is the sum of the  $k$-th column of $H(1,l,n)^{-1}$.     

It follows from (\ref{21.2}) that the terms in
(\ref{22.11j})--(\ref{22.19}) are all positive for all $l$ and $n$.    Since $H(1,l,n)^{-1}$ is also an M-matrix, it follows that $K^{-1}(1,l,n)  $ is a nonsingular M-matrix with all row sums equal to $0$,
except for the first row sum which is equal to 1. This completes the proof of Lemma \ref{lem-kmatrix}.
\qed

  We us the following notation.

\medskip	\noindent{\bf Notation }
 \begin{equation}
h_{j}= (f_{j}g_{j})^{1/2}, \hspace{.2 in}{\bf  v}(1,l,n)=(1,h_{2+l},\ldots, h_{n+l}),\label{ref-hV}
\end{equation}
\vspace{-.15 in}
\begin{equation}
m(1,l,n)_{j}= (r(1,l,n)_{j}c(1,l,n)_{j})^{1/2}, \,{\bf  m}(1,l,n)=(m(1,l,n)_{1},\ldots,m(1,l,n)_{n}),\label{ref-mM}
\end{equation}
\vspace{-.15 in}
 \begin{equation}
a(1,l,n)_{j} = \tau^{-1/2}_{1,l,n}\({\bf  m}(1,l,n) H(1,l,n)_{isymi}\)_{j}. \label{22.30hj}
 \end{equation}
 
 \medskip	For use below we note that $   r  (1,l,n)_{1}=1+o_{l}(1)$ and $c  (1,l,n)_{1}=   1+o_{l}(1)$ which imply that,
     \begin{equation}
  m(1,l,n)_{1}=1+o_{l}(1),\label{5.16}
    \end{equation}
and, 
\begin{equation}
   \sum_{j =2}^{n} m(1,l,n)_{j}\le   \(\sum_{j =2}^{n} r (1,l,n)_{j}\)^{1/2}\(\sum_{j =2}^{n} c (1,l,n)_{j}\)^{1/2}=o_{l}(1)\label{5.20mm}
   \end{equation}
 by (\ref{21.2}). 
 
   Note that for $0\le  a,b\le 1$, $(ab)^{1/2}\ge a\wedge b$, which gives $ (ab)^{1/2}\ge a+b-1$, or equivalently, $1-(ab)^{1/2}\le 2-a-b$. Consequently,
 \begin{equation}
    1-(f_{l+j}g_{l+j})^{1/2}\le  1- f_{l+j} + 1- g_{l+j}.
    \end{equation}
Therefore,   it also follows from (\ref{21.2}) that. 
    \begin{equation}
       1-h_{l+k}=o(\la_{l+k})\label{5.22mm}.
       \end{equation}

\medskip	\noindent{\bf  Proof of Lemma \ref{lem-koldet} } The lower bound is given in \cite[Lemma 3.3]{MRHD}.
We now obtain the upper bound. By subtracting the first line of (\ref{19.19e}) from each of the other lines and using (\ref{5.9mm}) we see that $|K(1,l,n)|=|H(1,l,n)|=\prod_{j =l+2 }^{l+n }\la_{j},     $  where we use (\ref{22.10}) for the last equality. Consequently,
  \begin{equation}
| \AA(1,l,n) |= | H(1,l,n)^{-1}|=\prod_{j =l+2 }^{l+n }1/\la_{j}.    \label{22.11}
 \end{equation}
 
Note that by and (\ref{2.15}) and (\ref{22.6}), 
     \be      H(1,l,n)_{sym}^{-1}=\left (
\begin{array}{ cccc  }  (1+\sum_{j =l+2 }^{l+n } h^{2}_{j}/\la_{j})&-h_{l+2}/\la_{l+2} &\dots&- h_{l+n }/\la_{l+n }  \\
-h_{l+2}/ \la_{l+2}&1/\la_{l+2}   & \dots &0   \\
\vdots&\vdots& \ddots&\vdots \\
-h_{l+n }/ \la_{l+n } &0 &  \dots &1/\la_{l+n }     \end{array}\right )\label{22.13}.
    \ee
     For $2\le j\le n$, multiply the $j-$th row of this matrix by $h_{l+j}$ and add it to the first row to see that,
    \be  | H(1,l,n)_{sym}^{-1}| =\prod_{j =l+2 }^{l+n }1/\la_{j}\label{22.14}.
 \ee

In our notation $\AA (1,l,n) =K(l,n)^{-1}$. Therefore, it follows from (\ref{22.18}) that,
    \be  \AA_{sym}(1,l,n)   =\left (
\begin{array}{ cccc  } 1+\rho_{1,l,n}&- m(1,l,n)_{1 }  &\dots&- m(1,l,n)_{n } \\
 - m(1,l,n)(l,n)_{1 }&H(1,l,n)_{sym}^{1 ,1 } & \dots &H(1,l,n)_{sym}^{ 1 ,n }   \\
\vdots&\vdots& \ddots&\vdots \\
  -m(1,l,n)_{n } &H(1,l,n)_{sym}^{n ,1 }&  \dots &H(1,l,n)_{sym} ^{n ,n }     \end{array}\right )\label{22.20}
     \ee 
We  write $\AA_{sym}(1,l,n) $ in block form
 \be \AA_{sym}(1,l,n) =\left (
\begin{array}{ ccccc }(1+\rho_{1,l,n}) &-{\bf  m}(1,l,n) \\
-{\bf  m}(1,l,n)^{T}&H(1,l,n)_{sym}^{-1}  \end{array}\right ),\label{22.21}
     \ee 
and use the formula for the determinant of a block matrix; see e.g.,   \cite[Appendix B]{DMM} to get,
 \bea
 &&
| \AA_{sym}(1,l,n) |\nn\\
&&\qquad=|H_{sym}(1,l,n)^{-1}|\,\((1+\rho_{1,l,n})-{\bf  m}(1,l,n)H(1,l,n)_{isymi}{\bf  m}(1,l,n)^{T}\),  \nn
 \eea
	and therefore by (\ref{22.11}) and (\ref{22.14}) we have,
 \be
\frac{ | \AA_{sym}(1,l,n) |}{  | \AA(1,l,n) |}   
= (1+\rho_{1,l,n})-{\bf  m}(1,l,n)H(1,l,n)_{isymi}{\bf  m}(1,l,n)^{T}.  \label{22.23} 
 \ee
By (\ref{con7qq})  
\bea 
   H(1,l,n)_{isymi} = H(1,l,n)_{Sym}=D_{\vec\la}+D_{\bf v}{\bf 1}D_{\bf v}.\label{5.26mm}
 \eea  
where $\vec \la=\(0,\la_{2+l},\ldots,\la_{n+l} \)$ and ${\bf v}={\bf v}(1,l,n)$. (See (\ref{1.35mm}) for notation.) 

To simplify the calculation we note that
\begin{equation}
   (D_{\bf v}{\bf 1}D_{\bf v})_{j+l ,k+l}= h_{j+l}h_{k+l},\qquad j,k=1,\ldots,n\label{5.27mm}
   \end{equation}
and $\la_{1+l}=0$ and  $h_{1+l}=1$. Consequently,
\bea
  && {\bf  m}(1,l,n)H(1,l,n)_{isymi}{\bf  m}(1,l,n)^{T}\\
  &&\qquad= \sum_{j,k=1}^{n} {\bf  m}(1,l,n)_{j} \(\la_{j+l}\de_{j+l,k+l}+h_{j+l}h_{k+l}\){\bf  m}(1,l,n)_{k}\nn\\
    &&\qquad= \sum_{j =2}^{n} \la_{j+l} m^{2}(1,l,n)_{j}       -\(m(1,l,n)_{1}+\sum_{j =2}^{n}  m(1,l,n)_{j}h_{j+l} \)^{2}\nn\\
&&\qquad\ge  \sum_{j =2}^{n} \la_{j+l} m^{2}(1,l,n)_{j}       -\ m^{2}(1,l,n)_{1} .\nn
   \eea
   Also, using (\ref{22.11j})-(\ref{22.19}), we have,
 \bea
 &&
 \sum_{j =2}^{n}\la_{j+l} m ^{2}(1,l,n)_{j}  =\sum_{j =2}^{n} \la_{j+l} r (1,l,n)_{j}  c(1,l,n) _{j} \label{22.26}\\
 &&\qquad = \sum_{j =2+l}^{n+l}\(1- g_{j} \)\(1- f_{j} \)\frac{1}{\la_{j}}=\rho_{1,l,n}-1.\nn
 \eea
Therefore,
\begin{equation}
   { | \AA_{sym}(1,l,n) | \over   | \AA(1,l,n) |}  \le 2-m^{2}(1,l,n)_{1} .
   \end{equation}
It follows from (\ref{5.16}) that both (\ref{22.15}) and (\ref{22.15a}) hold.\qed

The next lemma is used in the proof of Lemma \ref{lem-gold}
 
 \bl\label{lem-goldx} 
 Let $\eta_{n}$, $n=-1, 0,1,\ldots$ be a sequence of independent  standard normal random variables. The matrix $ \{K_{isymi}(1,l,n)_{j,k},j,k\in  1\ldots,n\}$ is the covariance matrix of $\{ \xi(1,l,n)_{j};j\in 1\ldots,n \}$ where,
 \begin{equation}
 \xi(1,l,n)_{1}=\eta_{-1}+ a(1,l,n)_{1}\eta_{0}, \label{22.321}
 \end{equation}
 and,
 \begin{equation}
\xi(1,l,n)_{j}= \la_{l+j}^{1/2} \eta_{l+j}+h_{l+j}\eta_{-1}+ a(1,l,n)_{j}\eta_{0}, \hspace{.2 in}j=2,\ldots,n, \label{22.32}
 \end{equation}
 for   a sequence  $ \{a(1,l,n)_{j}\}$ satisfying,
 \begin{equation}
a(1,l,n)_{j}=1+o_{l}(1)\qquad\mbox{and}\qquad a(1,l,n)_{j}-a(1,l,n)_{1}=o(\la_{j+l}),\label{22.50b}
\end{equation}
for $j=1,\ldots,n$.
  \el

\Proof  
To obtain (\ref{22.321}) and (\ref{22.32}) it suffices to show that,
\bea
K_{isymi}(1,l,n)_{1,1} &=&1+a(1,l,n)_{1}a(1,l,n)_{1}, \label{22.30aa1}\\&&\nn\\
K_{isymi}(1,l,n)_{1,k} &=&K_{isymi}(1,l,n)_{k,1} =h_{l+k} +a(1,l,n)_{1}a(1,l,n)_{k},\nn\\&&\nn\\
 K_{isymi}(1,l,n)_{j,k} &=&\la_{l+j}\de_{j,k}+h_{l+j}h_{l+k} +a(1,l,n)_{j}a(1,l,n)_{k}, \nn
 \eea
 $2\le j,k\le n.$
 
\medskip	 Since $K_{isymi}(1,l,n) =\( \AA_{sym}(1,l,n)\)^{-1}$  we use (\ref{22.21}) and the formula for the inverse of a matrix in  \cite[Appendix B]{DMM} to write  
  \be 
 K_{isymi}(1,l,n) =\left (
\begin{array}{ ccccc }C  &{\bf  w}  \\
{\bf  w}  ^{T}&\Phi(1,l,n) \end{array}\right ),\label{22.29}
     \ee
    for some constant $C$ and $n$-dimensional vector ${\bf  w} $. We do not need to know their precise values. What matters to us is the $n\times n$ matrix $\Phi (1,l,n)$ because it is $   \{K_{isymi}(1,l,n)_{j,k};j,k\in  1\ldots,n \}$ in (\ref{2.12nn}). We have,
 \bea
 \lefteqn{
\Phi (1,l,n)\label{22.30b}}\\
&&=H(1,l,n)_{isymi}+\tau^{-1}_{1,l,n} H(1,l,n)_{isymi}{\bf  m}(1,l,n)^{T}{\bf  m}(1,l,n) H(1,l,n)_{isymi}.\nn
 \eea
 
Let
 \begin{equation}
a(1,l,n)_{j} := \tau^{-1/2}_{1,l,n}\({\bf  m}(1,l,n) H(1,l,n)_{isymi}\)_{j}. \label{22.30}
 \end{equation}  
 Therefore, by (\ref{5.26mm}) and (\ref{5.27mm}),
 \begin{equation}
\Phi(1,l,n)_{j,k} = \la_{l+j}\de_{j,k}+  h_{l+j}  h_{l+k} +a(1,l,n)_{j}a(1,l,n)_{k}\label{22.30am}.
 \end{equation}
 Since $\la_{1+j}=0$ and $\la_{1+j}=0$, we get (\ref{22.30aa1}).

 \medskip	Using  (\ref{5.26mm}) and (\ref{5.27mm}) again we see that,
 
 \bea
  &&  ({\bf  m}(1,l,n) H(1,l,n)_{isymi})_{k}\\
   &&\qquad=\sum_{j=1}^{n}{\bf  m}(1,l,n)_{j} \(\la_{j+l}\de_{j+l,k+l}+h_{j+l }h_{ k+l}\)\nn\\
   &&\qquad= {\bf  m}(1,l,n)_{k}  \la_{k+l}+\sum_{j=1}^{n}{\bf  m}(1,l,n)_{j}  h_{j+l }h_{ k+l} \nn\\  &&\qquad= \(\la_{l+k}m(1,l,n)_{k}+\(m(1,l,n)_{1}+\sum_{j =2}^{n} m(1,l,n)_{j}h_{l+j} \)h_{l+k}\).\nn
    \eea
    Consequently, for all $k\in\mathbb N$, 
    \bea
&&a(1,l,n)_{k}\label{22.36}\\
&&\quad = \tau^{-1/2}_{1,l,n}\(\la_{l+k}m(1,l,n)_{k}+\(m(1,l,n)_{1}+\sum_{j =2}^{n} m(1,l,n)_{j}h_{l+j} \)h_{l+k}\).  \nn
\eea
In particular,
\begin{equation}
   a(1,l,n)_{1}= \tau^{-1/2}_{1,l,n} \(m(1,l,n)_{1}+\sum_{j =2}^{n}  m(1,l,n)_{j}h_{l+j} \) .
   \end{equation}
  It follows from (\ref{22.15a}), (\ref{5.16}), (\ref{5.20mm}) and the fact that $h_{j}\leq 1$ for all $j\in \mathbb N$ that 
 \begin{equation}
 a(1,l,n)_{k}=1+o_{l}\(1\),\qquad\forall j\in\mathbb N.\label{22.37}
 \end{equation}
  Using the last two displays we see that, for any $2\leq k\leq n$,   
\bea
&& a(1,l,n)_{k}-a(1,l,n)_{1}= \tau_{1,l,n}^{-1/2}\\
&& \times\(\la_{l+k}m(1,l,n)_{k}+\(m(1,l,n)_{1}+\sum_{j =2}^{n} m(1,l,n)_{j}h_{l+j} \)\( h_{l+k}-1\)\).\nn  
\eea
Therefore, by  (\ref{5.16}), (\ref{5.20mm}),  (\ref{5.22mm}) and the fact  $\tau_{1,l,n}^{-1/2}\leq 1$, we get,
\begin{equation}
a(1,l,n)_{k}-a(1,l,n)_{1}=o(\la_{k+l}).\label{22.50bv}
\end{equation}
 \qed

\noindent {\bf Proof of Lemma \ref{lem-gold} } It follows from Lemma \ref{lem-goldx} that for all $j\ge 1$,
\bea 
   a^{2}(1,l,n)_{j}= 1+o_{l}(1).
 \eea
Using this and (\ref{22.321}) we see that,
\bea
    \xi(1,l,n)_{1}&\stl& \(1+a^{2}(1,l,n)_{1}\)^{1/2}\rho\\
    &=&\sqrt 2\rho+o_{l}(1) \rho\nn,
   \eea
where $\rho$ is a standard normal random variable. Similarly, for $j=2,\ldots,n$,
\bea 
   h_{l+j}\eta_{-1}+ a(1,l,n)_{j}\eta_{0}&\stl& \(  h^{2}_{l+j}\eta_{-1}+a^{2}(1,l,n)_{j}\)^{1/2}\rho\\
    &=&\sqrt 2\rho+o_{l}(1) \rho\nn.
 \eea
Therefore, it follows from Lemma \ref{lem-goldx} that,
 \begin{equation}
\xi(1,l,n)_{j}\stl \la_{l+j}^{1/2} \eta_{l+j}+\sqrt 2\rho+o_{l}(1) \rho, \hspace{.2 in}j=2,\ldots,n, \label{22.32z}
 \end{equation}
 Furthermore, by (\ref{22.321})--(\ref{22.50b}) and (\ref{5.22mm}),
 \bea
  &&  \xi(1,l,n)_{j}- \xi(1,l,n)_{1}\\
  &&\qquad  \stl  \la_{l+j}^{1/2} \eta_{l+j}+(h_{l+j}-1)\eta_{-1}+ (a(1,l,n)_{j}-a(1,l,n)_{1})\eta_{0}\nn\\
  &&\qquad  \stl  \la_{l+j}^{1/2} \eta_{l+j}+\((h_{l+j}-1)^{2} + (a(1,l,n)_{j}-a(1,l,n)_{1})^{2}\)^{1/2} \rho\nn\\
  &&\qquad  \stl  \la_{l+j}^{1/2} \eta_{l+j}+ o( \la_{l+j})\rho\nn
    \eea
  Set $\rho=\eta_{0 }$, $\la_{l+j}=\la_{1/(l+j) }$ and $\eta_{l+j}=\eta_{1/(l+j) }$ and we
have Lemma \ref{lem-gold}.\qed

\section{Relation between the kernel $U$ and   Reuter's example }\label{sec-9}

 We explain how the examples of Markov chains introduced by Kolmogorov and Reuter led us to study permanental process with the kernel $U$ defined in (\ref{21.1c}). To begin  we  use  \cite[Theorem 4.1.3]{book} to define the  a class Markov chains by giving   their
$\al$-potentials.  The state space of the chains is the
sequence  $T = \{0, \frac{1}{2},
\frac{1}{3},\ldots, \frac{1}{n},
\ldots   \}$, with the topology inherited from the real line. Clearly $T$ is a
compact metric space with one limit point. Let $\{q_n\}^\infty_{n=2} $, $\{r_n\}^\infty_{n=2}$ and
$\{s_n\}^\infty_{n=2}$ be strictly positive real numbers satisfying,\be
\sum^\infty_{n=2} \frac{q_n}{r_n} < \infty, \quad \lim_{n\to
\infty} q_n = \infty, \quad  s_{n}<r_{n} \quad \mbox{ and}\quad \limsup_{n\to
\infty}   r_n-s_{n}<\ff.\quad
\label{q10.1}
\ee

  We   define an $\al$-potential $\{U^\alpha, \alpha>0\}$
on $C(T)$ in terms of its density $u^\alpha(x,y), x,y\in T$  with respect to
counting measure. Set,
\bea
   u^\alpha(0,0) &=& \frac{1}{\alpha +\alpha\sum\limits^\infty_{j=2}
\displaystyle     \frac{ s_j}{ \alpha+r_j}\frac{q_{j} }{r_j }},
\label{q10.3m}\\
 u^\alpha\left( 1/i, 0\right)&=&  u^\alpha
(0,0)\frac{r_i}{\alpha +r_i},\label{w10.3m}\\
     u^\alpha\left(0, 1/j\right) &=&    u^\alpha
(0,0)\frac{ s_j}{ \alpha+r_j}\frac{q_{j} }{r_j }
\label{r10.3m},\\
   u^\alpha\left(1/i, 1/j\right) &=& \delta_{ij}
\frac{1} { \alpha+r_j} + u^\alpha (0,0) \frac{r_i}{\alpha+r_i}
\frac{ s_j}{ \alpha+r_j}\frac{q_{j} }{r_j }. \label{x10.3m}
\eea

  The next two lemmas give properties of $U^\al$. They are proved in Section \ref{sec-app2a}.
\begin{lemma} \label{lem-7.1a} For $f$ bounded
\begin{equation} \label{}
   \lim_{i\to\ff}U^\alpha f (1/i)=U^\alpha f (0) .
\end{equation} 
In particular      $ U^\alpha
  $ takes $ C(T) \to C(T)$.  
\end{lemma}

\begin{lemma} \label{lem-9.1} The operator  
  $ U^\alpha$ satisfies, 
\bea
   \alpha U^\alpha 1 &=& 1\label{q10.5u},\\
    \|\alpha U^\alpha\|_{\ff}&\le &1,\label{q10.5v}\\
  \lim_{\alpha \to \infty} \alpha U^\alpha f(x)&=&  f(x), \qquad \forall x\in
T\label{q10.5ww}\\
   U^\alpha-U^\bb + (\alpha-\bb) U^\alpha U^\bb &=& 0,\label {q10.6} 
\eea 
  for $f\in C(T)$.
\end{lemma}

 \medskip  It
follows from \cite[Theorem 4.1.3]{book}
 that  $\{U^\alpha, \alpha>0\}$ is the
resolvent of a
Feller process $X$  with  state space
$T$. 
We now choose some $\al>0$ and consider   $\{u^\alpha(x,y), x,y\in T\}$ as the potential density of the Markov chain  $\ov X$   that is obtained by killing $X$ at an independent  exponential time with mean $1/\al$.

\begin{lemma} \label{}
    The $Q$ matrix for $\ov X$,      $\{Q_{x,y},x,y\in T \}$ is,
\be
\bordermatrix{&\phantom{q_2}&\phantom{q_2}&\phantom{q_2}&\phantom{q_2}
&\phantom{q_2}&  \cr
\hfill
  &-\infty&\displaystyle\frac{s_{ 2}q_2}{r_{ 2}} &\displaystyle\frac{s_{ 3}q_3}{r_{ 3}}&\displaystyle\frac{s_{ 4}q_4}{r_{ 4}}&\dots&  \cr
 &r_2&-(\al+r_2)&0&0 &\dots&\phantom{\cdot} \cr   &r_3&0&-(\al+r_3)&0 &\dots&  \cr   &r_4&0&0&-(\al+r_4)&\dots& \cr 
 &\vdots&\vdots&\vdots&\vdots&\ddots\cr}.\label{q10.14}
 \ee
\end{lemma}

\Proof   Let $\{L^y_t,t\in T \}$ be the local time of $\ov X$ at $y$, i.e., 
\begin{equation}
L^{y}_{t}=\int_{0}^{t}1_{\big\{\ov X_{s}=y\big\}}\,ds.\label{lt.1}
\end{equation}
Note that,
 \begin{equation}
u^{\al}(x,y)=E^{x}\(L^{y}_{\ff}\),\qquad \forall \,x,y\in T. \label{lt.2}
\end{equation}
Let 
\begin{equation}
A_{t}^{n}=\int_{0}^{t}1_{\{\ov X_{s}\in T_{n}\}}\,ds,\qquad \mbox{and}\qquad\ga_{t}^{n}=\sup \{s\,|\, A_{s}^{n}=t\},\label{lt.8}
\end{equation}
  and  define,  
\be
\ov X_{(n)}=\{\ov X_{(n), t};t\in T\}=\{\ov X_{\ga_{t}^{n}};t\in T\}.
\ee 
Note that $ \{\ov X_{(n), t};\,t\geq 0\}$ is the Markov chain on $T_{n}=\{0, {1}/{2},
 ,\ldots,  {1}/{n}\}$ obtained from $\ov X$ by deleting all times when $\ov X$ is in $\{T-T_{n}\}$. 
 Set
\begin{equation}
L^{y}_{n,t}=\int_{0}^{t}1_{\{\ov X_{(n),s}=y\}}\,ds,\quad\mbox{and}\quad u_{n}^{\al}(x,y)=E^{x}\(L^{y}_{n,\ff}\),\quad \forall\, x,y\in T_{n}.\label{lt.3}
\end{equation}
It follows from \cite[Corollary 1.32a]{F} that,
\begin{equation}
L^{y}_{n,t}=L^{y}_{\ga_{t}^{n}}, \quad \mbox{ and }\quad u_{n}^{\al}(x,y)=u^{\al}(x,y),\quad \forall x,y\in \,T_{n}.\label{lt.4}
\end{equation}

To find the $Q$ matrix of   $\ov X_{(n)}$   we first consider the  matrix 
  \be H(  n)  =\left (
  \begin{array}{ cccc  } u^{\al}(0,0) &  u^{\al}(0,1/2) &\dots& u^{\al}(0,1/n) \\
 u^{\al}(1/2,0) & u^{\al}(1/2,1/2) &\dots&  u^{\al}(1/2, 1/n)  \\
 \vdots&\vdots& \ddots&\vdots \\
 u^{\al}(1/n,0) & u^{\al}(1/n,1/2) &\dots&  u^{\al}(1/n,1/n) \end{array}\right ).\label{n22.5}\ee
  Set 
 \be 
 u=   u^\alpha(0,0),\quad f_{j}=\frac{r_j}{\alpha +r_j},\, \quad g_{j}=\frac{ s_j}{ \alpha+r_j }\frac{q_{j} }{r_j },\quad  \la_{j}=\frac{1}{\alpha +r_j}\label{n22.57}.
 \ee
 Then,   
      \be H(  n)  =\left (
\begin{array}{ cccc  } u & ug_{ 2 }  &\dots& ug_{ n } \\
 u f_{ 2 }&\la_{ 2 }+uf_{ 2 } g_{ 2 }& \dots &uf_{ 2 }g_{ n }  \\
\vdots&\vdots& \ddots&\vdots \\
u f_{ n } &uf_{ n }g_{ 2 }&  \dots &\la_{ n }+uf_{ n }g_{ n }   \end{array}\right ).\label{n22.5r}
     \ee 
      Subtract $f_{ j}$ times the first row of $H(n)$ from the $j$--th row of $H(n)$ for $j=2,\ldots, n$, to see that,
      \begin{equation}
| H( n) | =u\prod_{j = 2 }^{ n }\la_{j}>0.   \label{22.10s}
 \end{equation}
Therefore, $H( n)$ is invertible. One can check that,
      \be      H( n)^{-1}=\left (
\begin{array}{ cccc  }  (u^{-1}+\sum_{j = 2 }^{ n } f_{j}g_{j}/\la_{j})&-g_{ 2}/\la_{2 } &\dots&- g_{ n}/\la_{ n}  \\
 -f_{ 2}/ \la_{ 2}&1/\la_{ 2}   & \dots &0   \\
\vdots&\vdots& \ddots&\vdots \\
 -f_{ n}/ \la_{ n} &0 &  \dots &1/\la_{ n}     \end{array}\right ) \label{n22.6}
 \ee
 by computing $H( n)H( n)^{-1}$. 
 
Using (\ref{n22.57})  we see that, 
      \be  -Q_{n}:=     H( n)^{-1}=\left (
\begin{array}{ cccc  }  \(\displaystyle\frac{1}{u} +\sum_{j = 2 }^{n } \displaystyle\frac{q_{j}s_j}{ \alpha+r_j}\)&-\displaystyle\frac{s_{ 2}}{r_{ 2}}q_{ 2} &\dots&- \displaystyle\frac{s_{ n}}{r_{ n}}q_{ n}   \\
 -r_{ 2}&\alpha +r_{ 2}   & \dots &0   \\
\vdots&\vdots& \ddots&\vdots \\
 -r_{ n} &0 &  \dots &\alpha +r_{ n}     \end{array}\right ).\label{n22.6c}
    \ee
     Since  $Q_{n}$ is uniform and $-Q^{-1}_{n}=H( n)$, it follows that  $Q_{n}$ is the $Q$ matrix for $\ov X_{n}$. In addition,  since  $ \lim_{j\to\ff}{s_{j}}/{r_{j}} =1    $   by (\ref{q10.1}),
 \begin{equation} \label{}
   \frac{q_{j}s_j}{ \alpha+r_j}\sim    \frac{q_{j} }{ 1+\alpha/r_j},\qquad \mbox{as $\quad j\to\ff$.}
\end{equation}
 Since $\lim_{j\to\ff} q_j=\ff$, also by (\ref{q10.1}), we see that,  \begin{equation} \label{}
 \lim_{n\to\ff} \sum_{j = 2 }^{ n } \frac{q_{j}s_j}{ \alpha+r_j}=\ff. 
\end{equation}
This shows that 
$Q_{n}$ converges to the $Q$ matrix in (\ref{q10.14}). The fact that (\ref{q10.14}) is the 
$Q$ matrix for $\ov X$ then follows from \cite[Corollary 1.90 and Section 2.7]{F}.   \qed

\begin{remark} {\rm   One can also show  that (\ref{q10.14}) is 
  the $Q$ matrix for $\ov X$   without appealing to \cite{F}.   First of all, we see from (\ref{n22.6c}) that for any $1<j\le n $, $\ov X_{(n)}$ can only jump from $1/j$ to $0$. Letting $n\to\ff$ shows that  $\ov X$ can only jump from $1/j$ to $0$. Therefore,
 the holding time in $1/j$ and the probability of a jump from $1/j$ to $0$ are the same for $\ov X_{(n)}$ and $\ov X $.   Using (\ref{n22.57})  we get all the entries in  (\ref{q10.14}) except for the first row. 
 
 Let $q_{0,j}$ be the $j$--th entry in the first row of the $Q$ matrix for $\ov X$ and let   $H_{j}$,  $j>1$, be the time of the first jump of $\ov X $ from $0$ to $1/j$ . It follows from \cite[Chapter III, (57.14)]{RW} that
\begin{equation}
P^{0}\(L^{0}_{H_{j}}>t\,|\, H_{j}<\ff\)=\exp^{-t q_{0,j}} .\label{lt.6}
\end{equation}
  Let   $H_{n,j}$,  $j>1$, be the time of the first jump of $\ov X_{(n)} $ from $0$ to $1/j$ . Since  $\ga_{H_{n,j}}^{n}=H_{j}$ it   follows from   (\ref{lt.4}) that $L^{0}_{n,H_{n,j}}=L^{0}_{H_{j}}$.
Consequently, $q_{0,j}$ must be the same as the corresponding entry in (\ref{n22.6c}). This verifies all the entries in the first row of the $Q$ matrix for $\ov X$ except for the first entry.

To obtain the first entry of the $Q$ matrix for $\ov X$ we note that by (\ref{q10.1}),
\begin{equation} \label{}
  -q_{0,0}\geq \sum_{j=2}^{\ff}q_{0,j}=\sum_{j=2}^{\ff}\frac{s_jq_j}{r_{j}}=\ff,
\end{equation}
where we use the fact that $ \lim_{ j\to\ff}{s_{j}}/{r_{j}} = 1$.
 }
\end{remark}

\begin{remark} \label{}{\rm 
   The matrix in  (\ref{q10.14}) with $\al=0$ and $s_j=r_j$, $j=2,\ldots,\ff$, is the $Q$ of the Markov chain considered by Reuter in \cite{Reu}. Some additional background information is given in \cite[page 455]{book}}. 
\end{remark}

 We  define permanental sequences by kernels that are potential densities of Markov chains. Rather than use the potential densities in (\ref{q10.3m})--(\ref{x10.3m}) we can simplify things if 
 we write the  $\al$-potential $ U^\alpha $
on $C(T)$ in terms of a density   $\{v^\alpha(x,y), x,y\in T\}$  with respect to a finite measure $m$ on $T$ given by $m(1/n) = q_n/r_n$, $n\ge 2$, and
$m(0) =1$.
  That is, for all   $f\in C(T)$,
\be U^\alpha f(x) = \int\limits_T v^\alpha(x,y) f(y) m(dy)\label{q10.2},
\ee where
\bea
   v^\alpha(0,0) &=& \frac{1}{\alpha +\sum\limits^\infty_{j=2}
\displaystyle \frac{\alpha s_j}{ \alpha+r_j}\frac{q_{j}}{r_{j}}},
\label{q10.3}\\
 v^\alpha\left( 1/i, 0\right) &=&  v^\alpha
(0,0)\frac{r_i}{\alpha +r_i},\nn\\
     v^\alpha\left(0, 1/j\right)&=&    v^\alpha
(0,0)\frac{s_j}{\alpha+r_j}
\nn,\\
   v^\alpha\left(1/i, 1/j\right) &=& \delta_{ij}
\frac{1}{ \alpha+r_j}\frac{r_j}{q_j } + v^\alpha (0,0) \frac{r_i}{\alpha+r_i}
\frac{s_j}{\alpha+r_j}. \nn
\eea
Note that  $v^\alpha(x,y)$ is continuous at $0$.

\medskip
We make   further changes to simplify (\ref{q10.3}). 
We fix $\al=1/2$ and if necessary change the  $\{q_{j}\}$, $\{r_{j}\}$ and  $\{s_{j}\}$ so that  $v^{1/2}(0,0)=1$. That is, we require
\begin{equation}
 \sum\limits^\infty_{j=2}
\displaystyle \frac{  q_{j} }{r_j }\(\frac{   s_j}{ 1/2+r_j}\)=1.\label{q10.3s1}
\end{equation}
Next we set,
\begin{equation}
f_{i}=\frac{r_{i}}{\al+r_{i}}, \qquad g_{j}=\frac{s_{j}}{\al+r_{j}},\qquad \mbox{and}\qquad \la_{j}=\frac{r_{j}}{q_{j}(\al+r_{i})}\label{q10.3s}.
\end{equation}
  Let  $V \left(x, y\right)  =  v^{1/2}\left(x, y\right)$. With these changes      we can write (\ref{q10.3}) as,
\bea
V(0,0) &=& 1,
\label{q10.3v}\\
V\left( 1/i, 0\right) &=& f_{i},\nn\\
V\left(0, 1/j\right)&=& g_{j}
\nn,\\
V\left(1/i, 1/j\right) &=& \delta_{ij}
\la_{j} + f_{i}
g_{j}. \nn
\eea
It  follows from (\ref{q10.3s}) and  (\ref{q10.1}) that  for $j\geq 2$,
   \begin{equation}
   0<  f_{j}<1, \quad 0<  g_{j}<1, \mbox{ and }  \lim_{j\to\ff}f_{j}=1, \quad\lim_{j\to\ff}g_{j}=1, \quad\lim_{j\to\ff}\la_{j}=0.\label{q10.1j}
   \end{equation}
   We also have 
    \begin{equation}
    \sum ^\infty_{j=2}\frac{ (1-f_{j})}{\la_{j}}=\al  \sum ^\infty_{j=2} \frac{ q_{j}}{r_{j}}     <\ff,\label{q10.1j1}
       \end{equation}
   and 
    \begin{equation}
    \sum ^\infty_{j=2}\frac{ (1-g_{j})}{\la_{j}}=\al  \sum ^\infty_{j=2} \frac{ q_{j}}{r_{j}}(\al +r_{j}-s_{j})     <\ff.\label{q10.1j2}
       \end{equation}
       
  We now    choose values of $\{q_{j}\}$, $\{r_{j}\}$ and  $\{s_{j}\}$  so that in addition to  (\ref{q10.3s1}) we also have,
       \begin{equation}
       \sum ^\infty_{j=2}\frac{ (1-f_{j})}{\la_{j}}<1 \quad \mbox{and} \quad       \sum ^\infty_{j=2}\frac{ (1-g_{j})}{\la_{j}}<1.\label{q10.3s2}
       \end{equation}
       That is, we   choose values of $\{q_{j}\}$, $\{r_{j}\}$ and  $\{s_{j}\}$  so that,
          \begin{equation}
      \sum\limits^\infty_{j=2}
\displaystyle\frac{ q_{j}}{r_{j}}\( \frac{   s_j}{1/2+r_j}\)=1,\quad      \sum ^\infty_{j=2} \frac{ q_{j}}{r_{j}} <2 \quad \mbox{ and } \quad       \sum ^\infty_{j=2} \frac{ q_{j}}{r_{j}}(1/2 +r_{j}-s_{j})<2.\label{q10.3s3}
       \end{equation}
 We write,  \bea \label{9.69}
  &&\sum\limits^\infty_{j=2}
\displaystyle\frac{ q_{j}}{r_{j}}\( \frac{   s_j}{1/2+r_j}\)\\
&&\qquad = q\sum ^{j_0}_{j=2}\displaystyle\frac{ 1}{r_{j}}\( \frac{   s_j}{1/2+r_j}\)\nn + \sum _{j=j_0+1}^{\ff}\displaystyle\frac{ q_{j}}{r_{j}}\( \frac{   s_j}{1/2+r_j}\)
\eea  
 where we take    $r_j=s_j+\de$, $j\ge 2$,   for $\de<((r_j/4)+\(3/8\))\wedge 1$ and  set $ q_j=q $, for $j\in 2,\ldots,j_0$. We then choose $j_0$ such that the last sum in (\ref{9.69}) is less than $1 $ and then choose $q$ so that the sum  of the two terms is equal to 1.
 This gives the first statement in (\ref{q10.3s3}).
  To continue we note that it follows from  this equality that,
 \begin{equation} \label{}
  \sum\limits^\infty_{j=2}
\displaystyle\frac{ q_{j}}{r_{j}}\le \inf_{j\ge 2}\( \frac{   s_j}{1/2+r_j}\)^{-1}.
\end{equation}
  Therefore, if,
  \begin{equation} \label{9.71}
  \inf_{j\ge 2}\( \frac{   s_j}{1/2+r_j}\)>\frac{3}{4}
\end{equation}    
 we get 
 \begin{equation} \label{}
  \sum ^\infty_{j=2} \frac{ q_{j}}{r_{j}} <\frac{4}{3}. 
\end{equation} This gives the second relationship in (\ref{q10.3s3}). Clearly, the third  relationship  in (\ref{q10.3s3}) holds when in addition to the above,
 \begin{equation} \label{9.73}
 \sup_{j\ge 2}\( \frac{1}{2}+r_j-s_j\)\le \frac{3}{2}.
\end{equation}    
       Note that  (\ref{9.71}) and (\ref{9.73}) hold  for $\de<((r_j/4)+\(3/8\))\wedge 1$. 
 
\medskip   For $s,t\in T$,
  \begin{equation} \label{7.43}
  V\(s,t\)+\mathbf{1}(s,t)=U\(s,t\),
\end{equation}
  defined in  (\ref{21.1c}), where $\mathbf 1$ is an infinite matrix with all of its entries equal to 1. 
  
 \medskip  We now explain why 
  we don't simply work with the matrix $V$. 
  The reason is that although  $ V$ is clearly not symmetric,  it is  symmetrizable. 
   This is what this means.  Let $\mathbf i= i_{1},\ldots,i_{n}$, where $i_{1},\ldots,i_{n}\in T$. Let $D_{\mathbf { p_{  i} }}$ denote a diagonal matrix with diagonal elements $p_{i_{1}},\ldots,p_{i_{n}}$. Then   $\wt V_{\mathbf i}=  \{  V(i_{j},i_{k}), j,k=1,\ldots,n\}$ can be written as, 
\begin{equation}
  \wt V_{\mathbf i}=D_{\mathbf { \la_{  i} }}+D_{\mathbf { f_{  i} }}\mathbf 1D_{\mathbf {g_{  i} }}.\label{1.35mm}
   \end{equation}
Here $\mathbf { \la_{  i} }=\(\la (i_{1}),\ldots,  \la (i_{n})\)$ where      $ \la (0)=0, \,  \la (1/j)= \la_{j}$ for $j\geq 2$ and   $\mathbf { f_{  i} }=\(f(i_{1}),\ldots, f(i_{n})\)$  where      $f(0)=1,  \,  f(1/j)=f_{j}$ for $j\geq 2$. Similarly for $\mathbf {g_{  i} }$. Note that
 \begin{equation}
   |I+D_{\mathbf  s_{ i}}\wt V_{\mathbf i}|= |I+D_{s_{ \mathbf i}}  V^{*}_{\mathbf i}|,
    \end{equation}
for
 \begin{equation}
   V^{*}_{\mathbf i}=D_{\mathbf { \la_{  i} }}+D^{1/2}_{\mathbf { f_{  i} }}D^{1/2}_{\mathbf {g_{  i} }}\mathbf 1 D^{1/2}_{\mathbf { f_{  i} }}D^{1/2}_{\mathbf {g_{ i} }}.
   \end{equation}
  This is easy to see since, 
  \bea
       |I+D_{\mathbf  s_{ i}}\wt V_{\mathbf i}|&=& |I+D_{\mathbf  s_{ i}} (D_{\mathbf { \la_{  i} }}+D_{\mathbf { f_{  i} }}\mathbf 1D_{\mathbf {g_{  i} }})|\\&=&D^{1/2}_{\mathbf { f_{  i} }}   D^{-1/2}_{\mathbf {g_{  i} }}   |I+D_{\mathbf  s_{ i}} (D_{\mathbf { \la_{  i} }}+ D^{1/2}_{\mathbf { f_{  i} }}   D^{ 1/2}_{\mathbf {g_{  i} }}  \mathbf 1    D^{1/2}_{\mathbf { f_{  i} }}   D^{ 1/2}_{\mathbf {g_{  i} }} )|D^{- 1/2}_{\mathbf { f_{  i} }}   D^{ 1/2}_{\mathbf {g_{  i} }}\nn\\&=&   |I+D_{\mathbf  s_{ i}} (D_{\mathbf { \la_{  i} }}+ D^{1/2}_{\mathbf { f_{  i} }}   D^{ 1/2}_{\mathbf {g_{  i} }}  \mathbf 1    D^{1/2}_{\mathbf { f_{  i} }}   D^{ 1/2}_{\mathbf {g_{  i} }} )| = |I+D_{s_{ \mathbf i}}  V^{*}_{\mathbf i}|.
\nn
     \eea 
   
 Clearly, $V^{*}_{\mathbf i}$  is symmetric. We say that $  V$ is symmetrizable because  it determines the same $\al$-permanental sequences as a symmetric kernel. We also say that $  V$ is equivalent to $V^{*}_{\mathbf i}$. 
   
   Our goal is to study sample path properties of permanental  sequences that are not defined by symmetric kernels. Obviously   the permanental  sequence determined by $\{V(s,t);s,t\in T \}$ does not have this property.

 \medskip	   The kernel $U=\{U(s,t),s,t\in T\}$ defined in   (\ref{21.1c})   is generally not symmetrizable. 
 To see this it suffices to show that  the matrix $G$ defined in (\ref{1.2})
  is generally not  symetrizable. We actually want to show that for all   $l$,  the matrices,
  \begin{equation}
     G_{ l}:=\{G_{i,j};i,j\in (l,l+1,\ldots)\},
     \end{equation}
  are not symetrizable. This is because  for the asymptotic  results we obtain  are     only concerned with the   kernel $U(s,t)$ as $s,t\to 0$. 
 
  Note that if $G_{ l}$ is symmetrizable then  the $3\times 3$ matrices $G(3 l,3)$, the matrices  $G$ restricted  to
 $\{3l+1,3l+2,3l+3\}\times \{3l+1,3l+2,3l+3\}$, are symmetrizable for all $l$.
 We now show that when the $\{g_{i}\}$ are all different   we can find  $0<f_{j} <1$ such that,
\begin{equation}
0<(1- f_{j})< {\ep  \la_{j}\over 2^{j+1}},\qquad\forall  j\geq 2,   \label{nos.1}
\end{equation} 
  for any $\ep>0$, with the property that   for all $l\geq 1$ the $3\times 3$ matrix $G(3 l,3)$  is not equivalent to a symmetric matrix. It follows from   \cite[Lemma 2.1, $(2.6)$]{MRns} that for $G(3 l,3)$ to be equivalent to a symmetric matrix,   we must have  
 \be  ( 1+f_{i_{1}}g_{i_{2}})( 1 +f_{i_{2}}g_{i_{3}})( 1+f_{i_{3}}g_{i_{1}})  =( 1+f_{i_{1}}g_{i_{3}})( 1+f_{i_{2}}g_{i_{1}})( 1+f_{i_{3}}g_{i_{2}}),  \label{3.9wq}
  \ee  
where   $(i_{1},i_{2},i_{3})=(3l+1,3l+2,3l+3)$.
  Assume that this holds for all   values of $0<f_{i_{1}},  f_{i_{2}}, f_{i_{3}}<1$ satisfying (\ref{nos.1}). 
Think of   $f_{i_{1}}, f_{i_{2}}$ and $ f_{i_{3}}$  as real variables. If (\ref{3.9wq}) holds then its derivatives with respect to $f_{i_{1}}$ and $f_{i_{2}}$ must be equal. That is, we must  have, 
\be
g_{i_{2}}g_{i_{3}}( 1 +f_{i_{3}}g_{i_{1}})= g_{i_{1}}g_{i_{3}}( 1 +f_{i_{3}}g_{i_{2}}).\label{6.12}
\ee 
For (\ref{6.12}) to hold we must have $ g_{i_{2}}=g_{i_{1}}$. However, we start with all $\{g_{i}\}$  different. Consequently, (\ref{3.9wq})  does not hold for all  $0<f_{i_{1}},  f_{i_{2}}, f_{i_{3}}<1$ satisfying (\ref{nos.1}).   Therefore, we can choose   $0<f_{3l+1},  f_{3l+2}, f_{3l+3}<1$
satisfying (\ref{nos.1}), that do not satisfy (\ref{3.9wq}). Obviously we can choose 
$ f_{3l+1},  f_{3l+2}, f_{3l+3} $ so that in addition, $f_{j}\neq g_{j}$.  

\medskip Finally we have the following important observation:  
 
\begin{theorem} \label{theo-8.1} The matrix
  $\{U(s,t),s,t\in T\}$ is the restriction to $T\times T$ of  the potential density of a transient Markov chain with state space $\wt T=T\cup \{\ast\}$, where $\ast$ is  an isolated point.  
  \end{theorem}

To prove this theorem we use the following lemma: 
 
 \begin{lemma} \label{lem-7.4} The vector $\vec 1=(1,1,\ldots)$ is a left--potential for $\{V(x,y);x,y\in T \}$. That is there exists a function $h(x), x\in T$ with  $  h>0$ and $1<\|  h\|_1<2$, such that,
     \begin{equation}
 \sum_{x\in T}  h(x) V(x,y)=1,\qquad \forall y\in T.\label{n11.31a}
     \end{equation}

   \end{lemma}
   
 \Proof  Let
\be   
h(0)=1-   \sum_{j = 2 }^{\ff}f_{j} \( 1-g_{j}\)\frac{1}{\la_{j}},\label{n22.13j}\quad \mbox{and}\quad  h(1/k) =\(1-g_{ k} \)\frac{1}{\la_{ k}}.\quad k\geq 2.
\ee
By (\ref{q10.1j}) and (\ref{q10.3s2}) we have that  $h(x)> 0$ for all $x\in T$ and,  
\bea
\sum_{x\in T }^{\ff} h( x)& =&h(0) +  \sum_{j = 2 }^{\ff} h(1/j)\\
&=&1-   \sum_{j = 2 }^{\ff}f_{j} \( 1-g_{j}\)\frac{1}{\la_{j}}  +  \sum_{j = 2 }^{\ff}  \( 1-g_{j}\)\frac{1}{\la_{j}}<2.\nn
\eea
 This shows that $\|  h\|_1<2$.
 
 Now refer to (\ref{n22.13j}) and note that,\begin{equation} \label{}
  h(0)=1-   \sum_{j = 2 }^{\ff}f_{j}h(1/j). 
\end{equation}  
Using this and (\ref{q10.3v}) we see that,
  \begin{equation}
   \sum_{x\in T}h(x) V(x,0)= h(0) +  \sum_{j = 2 }^{\ff}f_{j} h(1/j) =1.\label{n11.40}
     \end{equation}
    In addition for $y=1/k$, $k\geq 2$,
    \begin{eqnarray}
     \sum_{x\in T}h(x) V(x,1/k)&=&h(0)V(0,1/k)+ \sum_{j= 2}^{\ff}h(1/j) V(i/j,1/k) \label{n11.41}\\
    &=&h(0)g_{k} + h(1/k)\la_{k}+ \sum_{j = 2 }^{\ff} h(1/j) f_{j}g_{k}\nn
   \\
    &=&g_k\(h (0)+\sum_{j = 2 }^{\ff}f_{j} h(1/j) \)+h(1/k)\la_{k} \nonumber\\
     &=&g_k +h(1/k)\la_{k}=1\nn,
    \end{eqnarray} 
  which gives (\ref{n11.31a}).
     
   Let $f(0)=1$ then since,
 \begin{equation} \label{}
  \sum_{x\in T}f(x)h(x)=1,
\end{equation}
and $f(x)<1$ when $x\ne 0$, $\|h\|_{1}>1$.
  \qed

 \medskip\noindent {\bf Proof of Theorem \ref{lem-7.4} }  The matrix $\{V(s,t);s,t\in T \}$ is the potential of a transient Markov chain on $T$. It follows from Lemma \ref{lem-7.4} that $\vec 1$ is a left-potential for this chain.
 Therefore, by \cite[Theorem 6.1]{MRsp}    there exists a transient Markov chain 
 $\wt X$ with state space $\wt T=T\cup \{\ast\}$, where $\ast$ is  an isolated point, such that $\wt X$ has potential densities,  
\be
   \wt    V(s,t)= V(s,t)+1=U(s,t), \hspace{.2 in}s,t\in T,
\label{rp.2} 
\ee
and
\be   \wt    V(\ast,t)= \wt    V(s, \ast)=\wt    V(\ast, \ast) =1, 
\nonumber
\ee 
with respect to a finite measure $\wt m$ on $\wt T$ which is equal to $m$ on $ T$ and has  $\wt m (\ast)=1$.   \qed

\begin{remark}  {\rm
It   follows from Theorem \ref{theo-8.1} and \cite[Theorem 3.1]{EK} that $\wt    V=\{\wt    V(s,t),s,t\in \wt T\}$ is the kernel of an $\al$ permanental process for all $\al>0$.     Theorem \ref{theo-1.1mm} then follows from the fact that    $U=\{U(s,t),s,t\in T\}$ is the restriction of $\wt    V=\{\wt    V(s,t),s,t\in \wt T\}$ to   $\{s,t\in T\}$. (However, the proof of Theorem \ref{theo-1.1mm} in Section 5 is much simpler than this.)}
  
\end{remark}

\section{Appendix I}\label{sec-app1}	

\begin{lemma} \label{lem-7.1}
Let $H=\{H_{j,k};j,k=1,\ldots,n\}$ be an $n\times n$    nonsingular matrix and set, 
\bea 
   G_{j,k}=H_{j,k}+h_{k},\qquad j,k=1,\ldots,n,
 \eea
  and,
  \be K  =\left (
\begin{array}{ ccccc }  1 & h_{1}&\ldots& h_{n}  \\
1   &G_{1,1} &\ldots&G_{1,n}   \\
\vdots& \vdots &\ddots &\vdots  \\
1   &G_{n,1} &\ldots&G_{n,n} 
\end{array}\right ). \label{19.19ewq}
     \ee  
Then $K$ is nonsingular and,
 \be K^{-1}  =\left (
\begin{array}{ ccccc }  1+\sum_{j,k=1}^{n}h_{k} H^{j,k} &-\sum_{j=1}^{n}h_{j}H^{j,1}&\ldots&-\sum_{j=1}^{n}h_{j}H^{j,n}  \\&&\vspace{-.1 in}\\
-\sum_{k=1}^{n} H^{1,k}   &H^{1,1} &\ldots&H^{1,n}   \\
\vdots& \vdots &\ddots &\vdots  \\
-\sum_{k=1}^{n} H^{n,k}   &H^{n,1} &\ldots&H^{n,n} 
\end{array}\right ). \label{19.19ew}
     \ee  
     
 Furthermore,    if  $H^{-1}$ is an $M$-matrix with positive row sums and 
 \begin{equation}
    \sum_{j=1}^{n}h_{j}H^{j,k}\ge 0,\qquad \forall k\in 1,\ldots,n,
    \end{equation}
    then $K$ is an inverse $M$-matrix. 
     \end{lemma}
     
     Note that all the row sums of $K^{-1}$ are equal to 0, except for the first row sum which is equal to 1.

 \medskip	\Proof Subtract  the first row of $K$ from each of the other rows to see that $|K|=|H|$. Therefore, $K$ is nonsingular. To obtain (\ref{19.19ew}) we   show that
   $ABCK=I$, where 
  	
 \be C  =\left (
\begin{array}{ ccccc }  1 & 0&0&\ldots& 0  \\-1   &1 & 0&\cdots&0   \\
-1   &0 & 1&\ldots&0   \\
\vdots& \vdots& \vdots  &\ddots &\vdots  \\
-1   &0&0 &\ldots&1
\end{array}\right ) \label{19.19ewa} .
     \ee  
and $A$ and $B$ and $C$, written in block form are,
 \be 
 A =\left (
\begin{array}{ ccccc }1  &{-\bf  h} \\
{\bf  0^{T}} &I \end{array}\right ), \qquad B =\left (
\begin{array}{ ccccc }1  &{\bf  0} \\
{\bf  0^{T}} &H^{-1} \end{array}\right )\label{22.29q},\qquad C =\left (
\begin{array}{ ccccc }1  &{ \bf  0} \\
({-\bf  1})^{T} &I \end{array}\right ) \qquad
     \ee
where ${\bf  h} $,   ${\bf  0} $, and ${\bf  1} $ are $n$-dimensional vectors.
  
 			It is easy to see that $ABCK=I$. The operation $CK$
subtracts the first row of $K$ from each of the other rows. Therefore
 \be 
 CK =\left (
\begin{array}{ ccccc }1  &{\bf  h} \\
{\bf  0^{T}} &H  \end{array}\right ).\label{22.29qw}
     \ee
Then
 \be 
 BCK =\left (
\begin{array}{ ccccc }1  &{\bf  h} \\
{\bf  0^{T}} &{I }  \end{array}\right ),\label{22.29qwd}
     \ee
where $I$ is an $n\times n$ identity matrix. Applying $A$ to this replaces $\bf h$ by $\bf 0$.

 Therefore, $K^{-1}=ABC$. Note that 
\begin{equation}
B C=\left (
\begin{array}{ ccccc }1  &{\bf  0} \\
   - H^{-1}  {\bf 1}^{T} &H^{-1} \end{array}\right )
   \end{equation}
and
\begin{equation}
AB C=\left (
\begin{array}{ ccccc }1+{\bf h}\cdot     (H^{-1}) {\bf 1 }^{T} &-{\bf h}H^{-1} \\
    -(H^{-1}) \bf 1 ^{T} & H^{-1} \end{array}\right ).
   \end{equation}
This is (\ref{19.19ew}).\qed

 \section{Appendix II}\label{sec-app2}
 
 \noindent{\bf Proof of Lemma \ref{lem-7.1a} }  
We have,
 \bea
 \lefteqn{ U^\alpha f (1/i)=\sum_{x\in T}u^\al
(1/i,x)f(x)=u^\al
(1/i,0)f(0)+\sum_{j=2}^{\ff }u^\al
(1/i,1/j)f(1/j)}\\&&=u^\alpha
(0,0)\frac{r_i}{\alpha +r_i}f(0)+\frac{1} { \alpha+r_i}f (1/i) +\sum_{j=2}^{\ff} u^\alpha (0,0) \frac{r_i}{\alpha+r_i}
\frac{ s_j}{ \alpha+r_j}\frac{q_{j} }{r_j }f(1/j)\nn .
\eea
Therefore, using  (\ref{q10.1}) and the Dominated Convergence Theorem we see that, 
\begin{equation} \label{}
  \lim_{i\to\ff}U^\alpha f (1/i)= u^\alpha
(0,0) f(0)+\sum_{j=2}^{\ff} u^\alpha (0,0)  
\frac{ s_j}{ \alpha+r_j}\frac{q_{j} }{r_j }f(1/j).
\end{equation}
The reader can check that this is equal to $U^\al f(0)$. The second statement in the lemma follows immediately from the first one.\qed

\noindent{\bf Proof of Lemma \ref{lem-9.1}\, }\label{sec-app2a} Recall that,
\begin{equation} \label{}
  U^{\al}f(x)=\sum_{x\in T}u^{\al}(x,y)f(y).
\end{equation}
We show below that for all $x\in T$,
\begin{equation} \label{9.7}
  \sum_{y\in T}u^{\al}(x,y)=\frac{1}{\al}.
\end{equation}
This immediately gives (\ref{q10.5u}) and (\ref{q10.5v}) since, 
\begin{equation} \label{}
   \sum_{y\in T}u^{\al}(x,y)f(y)\le  \sum_{y\in T}u^{\al}(x,y)\|f \|_{\ff}. 
\end{equation}

\medskip To obtain (\ref{9.7}) we note that by (\ref{q10.3m}),
\begin{equation} \label{}
  \frac{1}{\al u^{\al}(0,0)}=1+\frac{1}{  u^{\al}(0,0)} \sum_{j=2}^{\ff}u^{\al}(0,1/j).
\end{equation}
 Multiplying by $u^{\al}(0,0)$
  gives (\ref{9.7}) when $x=0$.
In addition, using  (\ref{r10.3m}) and (\ref{x10.3m}) we see that,
\be  \label{}
  \sum_{j=2}^{\ff}u^{\al}(1/i,1/j) =\frac{1}{\al+r_i}+ \frac{r_{i} }{\al+r_i}  \sum_{j=2}^{\ff}u^{\al}(0,1/j).    
   \ee 
Consequently adding $u^{\al}(1/i,0)$ to each side and using (\ref{w10.3m}) and (\ref{9.7}) in the case $x=0$ we get,  
\be  \label{9.16}
  \sum_{y\in T} u^{\al}(1/i,y) =\frac{1}{\al+r_i}+ \frac{r_{i} }{\al+r_i}  \sum_{y\in T} u^{\al}(0,y)=\frac{1}{\al},    
   \ee 
which gives (\ref{9.7})   for the other values of  $x$.

\medskip To obtain  (\ref{q10.5ww}) we first note that by the Dominated Convergence Theorem,
 \begin{equation} \label{}
  \lim_{\al\to\ff} \sum\limits^\infty_{j=2}
    \displaystyle  \frac{ s_j}{ \alpha+r_j} \frac{q_{j} }{r_j }=0.
\end{equation}
Therefore, it follows from (\ref{q10.3m}) that,
 \begin{equation}
\lim_{\al\to\ff}\al u^{\al}(0,0)=  1. \label{35.1}
 \end{equation} 
 Using this and the Dominated Convergence Theorem again, we see that  for any $f\in C\(T\)$,
\begin{eqnarray}
   \lim_{\al\to\ff}  \al U^\alpha f(0)  &=& \lim_{\al\to\ff}\al  u^{\al}(0,0)f(0)+ \lim_{\al\to\ff}\al\sum_{j=2}^{\ff}u^{\al}(0,1/j)f(1/j)
\nn\\
&=& f(0)+  \lim_{\al\to\ff} \sum_{j=2}^{\ff}\frac{s_j}{\al+r_j }\frac{q_j }{r_j}f(1/j)=f(0).\nonumber
\end{eqnarray}
  Similarly, for $i\ge 2$,    
\be  \label{9.33}
   \lim_{\al\to\ff} \al U^\alpha f(1/i) =\lim_{\al\to\ff}\al u^{\al}(1/i,0)f(0)+\lim_{\al\to\ff}\al\sum_{j=2}^{\ff}u^{\al}(1/i,1/j)f(1/j).  
  \ee 
Using (\ref{35.1}) we see that,  
  \begin{equation}
 \lim_{\al\to\ff}\al u^{\al}(1/i,0)f(0)=\lim_{\al\to\ff} \al u^{\al}(0,0)\frac{ r_i}{ \alpha+r_i}f(0)=0, \label{}
  \end{equation}
and,
\begin{eqnarray}
 &&\lim_{\al\to\ff}\al\sum_{j=2}^{\ff}u^{\al}(1/i,1/j)f(1/j) \\
&& \qquad =\lim_{\al\to\ff}\frac{ \al}{ \alpha+r_i}f(1/i)  +\lim_{\al\to\ff} \sum_{j=2}^{\ff}   \frac{ r_i}{ \alpha+r_i} \frac{s_i}{\al+r_i}\frac{q_i}{r_i}  f(1/j) = f(1/i)\nonumber.
\end{eqnarray}

\medskip  To obtain (\ref{q10.6}) it suffices to show that,
\be u^\al(x,y)-u^\bb(x,y)=(\bb-\al)\sum_{z\in T} u^\alpha(x,z) u^\bb(z,y),
 \label{q10.5f}
\ee for all $x,y\in T$. 
We first show this when $x=y=0$. We have,
\begin{eqnarray}
  \frac{u^\al(0,0)-u^\bb(0,0)}{u^\al(0,0)u^\bb(0,0)}&=& \(\bb\(1+\sum\limits^\infty_{j=2}
\displaystyle\frac{ s_j}{ \bb+r_j}\frac{q_{j} }{r_j }\)- \alpha\(1+\sum\limits^\infty_{j=2}
\displaystyle\frac{ s_j}{ \alpha+r_j}\frac{q_{j} }{r_j } \)  \) \label{res.12} \nn\\
 &=&\bb-\al+
    \sum\limits^\infty_{j=2}
\displaystyle\frac{  q_{j}s_j}{r_j }\(\frac{\bb}{(\bb+r_j)}-\frac{\al}{(\al+r_j)}\)\nn\\
 &=&(\bb-\al) \(1+    \sum\limits^\infty_{j=2}\frac{ s_jq_{j}}{ (\alpha+r_j)(\bb +r_j)}\) .\nonumber 
\eea
This last term,
\be  =(\bb-\al) \(1+   \frac{ 1}{u^\al(0,0)u^\bb(0,0)}\sum\limits^\infty_{j=2}u^\al(0,1/j)   u^\bb(1/j,0) \) .
\ee
Multiplying by $u^\al(0,0)u^\bb(0,0)$ we get (\ref{q10.5f}) when $x=y=0$.

\medskip We now obtain (\ref{q10.5f}) when $x=1/i$ and $ y=1/j$, $i,,j=2,\ldots.$
     By (\ref{r10.3m}) and (\ref{x10.3m}),
\begin{equation} \label{marc.1}
   u^{\al}(1/i,1/k)=\frac{\de_{i,k}}{\al+r_{k}}+ \frac{u^{\al}(1/i,0)}{u^{\al}(0,0)}u^{\al}(0,1/k),
\end{equation}
and,
\begin{equation} \label{marc.2}
   u^{\bb}(1/k,1/j)=\de_{j,k}\frac{1}{\bb+r_{j}}+\frac{u^{\bb}(1/k,0)}{u^{\bb}(0,0)}u^{\bb}(0,1/j).
\end{equation} 
Therefore, when $  i\ne j$,  
\bea \label{}
    &&u^{\al}(1/i,1/j)-u^{\bb}(1/i,1/j)\\  &&\quad =\frac{\de_{i,j}}{\al+r_{j}}- \frac{\de_{i,j}}{\bb+r_{j}}+\frac{u^{\al}(1/i,0)}{u^{\al}(0,0)}u^{\al}(0,1/j)-\frac{u^{\bb}(1/i,0)}{u^{\bb}(0,0)}u^{\bb}(0,1/j)\nn
\eea\vspace{-.2in}
\bea \label{}
 &&\nn = \frac{(\bb-\al)\de_{i,j}}{(\al+r_{j})(\bb+r_{j})} +\(\frac{1}{u^{\bb}(0,0)}-\frac{1}{u^{\al}(0,0)}\) u^{\al}(1/i,0) u^\bb(0,1/j)\\
 &&\qquad  \nn+ u^{\al}(1/i,0)\(\frac{u^\al(0,1/j)}{u^\al(0,0)}-\frac{u^\bb(0,1/j)}{u^\bb(0,0)}\)\\
  &&\qquad \quad \nn+ u^{\bb}( 0,1/j)\(\frac{u^\al( 1/i,0)}{u^\al(0,0)}-\frac{u^\bb( 1/i,0)}{u^\bb(0,0)}\):=a+b+c+d.
\eea
In addition,
\bea \label{}
  &&\sum_{z\in T}u^{\al}(1/i,z)u^{\bb}(z ,1/j)\\
  &&\qquad =u^{\al}(1/i,0)u^{\bb}(0,1/j)+  \sum_{k=2}^\ff u^{\al}(1/i,1/k)u^{\bb}(1/k,1/j).\nn
\eea 
We now multiply (\ref{marc.1}) and (\ref{marc.2}) and sum to obtain four terms in
 \begin{equation} \label{}
  \sum_{k=2}^{\ff} u^{\al}(1/i,1/k) u^{\bb}(1/k,1/j) .
\end{equation}
We have,  
\begin{equation} \label{}
  \sum_{k=2}^{\ff}\frac{\de_{i,k}}{\al+r_{k}}\frac{\de_{j,k}}{ \bb+r_{k}}=\frac{\de_{i,j}}{(\al+r_{i})(\bb+r_{j})}:=i
\end{equation}
\begin{equation} \label{}
  \sum_{k=2}^{\ff} \frac{\de_{i,k}}{\al+r_{k}}\frac{u^{\bb}(1/k,0)}{u^{\bb}(0,0)}u^{\bb}(0,1/j)=  \frac{ u^{\bb}(1/i,0)u^{\bb}(0,1/j)}{(\al+r_{i})u^{\bb}(0,0)}:= ii
\end{equation}
\begin{equation} \label{}
  \sum_{k=2}^{\ff} \frac{\de_{j,k}}{\bb+r_{j}}\frac{u^{\al}(1/i,0)}{u^{\al}(0,0)}u^{\al}(0,1/k)= \frac{ u^{\al}(1/i,0)u^{\al}(0,1/j)}{(\bb+r_{j})u^{\al}(0,0) }:=iii
\end{equation}
and since,
\bea \label{}
  &&\sum_{k=2}^{\ff}u^{\al}(0,1/k) u^{\bb}(1/k,0)=
\sum_{z\in T}^{\ff}u^{\al}(0,z) u^{\bb}(z,0)-u^{\al}(0,0) u^{\bb}(0,0)\nn\\
&&\qquad= \frac{u^{\al}(0,0)-u^{\bb}(0,0)}{\(\bb-\al\)}-u^{\al}(0,0) u^{\bb}(0,0),
\eea
we have,  
\bea \label{}
 &&\frac{u^{\al}(1/i,0)}{u^{\al}(0,0)}\frac{u^{\bb}(0,1/j)}{u^{\bb}(0,0)} \sum_{k=2}^{\ff}u^{\al}(0,1/k) u^{\bb}(1/k,0)\\
 &&\qquad  =u^{\al}(1/i,0)u^{\bb}( 0,1/j)\( \frac{u^{\al}(0,0)-u^{\bb}(0,0)}{\(\bb-\al\)u^{\al}(0,0)u^{\bb}(0,0)}-1\)\nn
  \eea
  We add $u^{\al}(1/i,0)u^{\bb}( 0,1/j)$ to this to get,
\bea \label{}
  \sum_{z\in T}u^{\al}(1/i,z)u^{\bb}(z ,1/j)&=& \frac{u^{\al}(1/i,0)u^{\bb}( 0,1/j)}{\(\bb-\al\)}\( \frac{1}{ u^{\bb}(0,0)}- \frac{1}{u^{\al}(0,0) }\)\nn\\
  & &\hspace{-.15in}+i+ii+iii:=iv+i+ii+iii.
\eea 
Note that,
\begin{equation} \label{}
  (\bb-\al)iv=b,\quad (\bb-\al) i =a
\end{equation}
 Also,
 \bea \label{}
  c&=&u^{\al}(1/i,0)\frac{q_j}{r_j}\(\frac{s_j}{\al+r_j}-\frac{s_j}{\bb+r_j}\)\\
  &=& (\bb-\al)u^{\al}(1/i,0) \(\frac{s_j }{(\al+r_j)}\frac{q_j}{r_j}\)\(\frac{1}{ \bb+r_j } \)=(\bb-\al)iii.\nn
\eea
Similarly, $(\bb-\al)ii=d$. This completes the proof of (\ref{q10.6}) for $x=1/i$ and $y=1/j$. The cases where $x=0$ and $y=1/j$, and $x=1/i$ and $y=0$ are similar.\qed

\noindent
\begin{tabular}{lll} & \hskip20pt Michael B.  Marcus
     & \hskip20pt  Jay Rosen\\  & \hskip20pt 253 West 73rd. St., Apt. 2E
   & \hskip20pt Department of Mathematics \\    &\hskip20pt  New York, NY 10023, USA
 \hskip20pt 
     & \hskip20pt College of Staten Island, CUNY \\    & \hskip20pt mbmarcus@optonline.net
    & \hskip20pt  Staten Island, NY
10314, USA \\    & \hskip20pt      & \hskip20pt   jrosen30@optimum.net \\   & \hskip20pt 

     & \hskip20pt 
\end{tabular}

  \end{document}